\documentclass[12pt]{article}  
\usepackage{graphicx,amsmath,amsthm,amssymb,dsfont}  
\numberwithin{equation}{section}
\def\RR{\mathbb R}
\def\CC{\mathbb C}
\def\ZZ{\mathbb Z}
\def\PP{\mathbb P}
\def\MM{\mathbb M}
\def\QQ{\mathbb Q}

\def\FF{\mathbb F}

\def\cK{\mathcal K}
\def\cR{\mathcal R}

\def\cD{\mathcal D}
\def\cQ{\mathcal Q}
\def\cC{\mathcal C}
\def\cA{\mathcal A}

\def\cL{\mathcal L}

\def\cM{\mathcal M}
\def\cJ{\mathcal J}
\def\cI{\mathcal I}
\def\cB{\mathcal B}

\def\cF{\mathcal F}

\def\fa{\mathfrak a}
\def\fm{\mathfrak m}

\def\cE{\mathcal E}

\def\fO{\mathcal O}
\def\cP{\mathcal P}
\def\fp{\mathfrak p}

\def\cG{\mathcal G}
\def\cD{\mathcal D}

\def\cN{\mathcal N}
\def\fN{\mathfrak N}

\def\cU{\mathcal U}
\def\fU{\mathfrak U}

\def\cS{\mathcal S}
\def\cY{\mathcal Y}
\def\cX{\mathcal X}
\def\cW{\mathcal W}

\def\bE{{\bf E}}

\def\bt{{\bf t}}
\def\bx{{\bf x}}
\def\by{{\bf y}}
\def\bz{{\bf z}}

\def\bk{{\bf k}}
\def\bl{{\bf l}}
\def\bp{{\bf p}}

\def\ba{{\bf a}}
\def\bb{{\bf b}}
\def\bq{{\bf q}}
\def\bff{{\bf f}}
\def\be{{\bf e}}
\def\bN{{\bf N}}
\def\Nm{{\rm Nm}}
\def\bgamma{\boldsymbol{\gamma}}
\def\btheta{\boldsymbol{\theta}}

\def\bvarsigma{\boldsymbol{\varsigma}}
\def\bepsilon{\boldsymbol{\epsilon}}
\def\bs{{\bf 0}}
\def\b1{{\bf 1}}
\def\d1{\mathds{ 1}} 
\def\sgn{{\rm sgn}}
\def\vol{{\rm vol}}

\def\qed{ \ \vrule width.2cm height.2cm depth0cm\smallskip}

\begin{document}
\hsize=14.8 true cm



\title{On the lower bound in the  lattice point remainder problem  for a parallelepiped}
\author{Mordechay B. Levin}

\date{}

\maketitle

\begin{abstract}
Let $ \Gamma \subset \RR^s $ be a lattice,  obtained from a 
 module in a totally real algebraic number field.  Let $G$ be an axis parallel  parallelepiped, and let $|G|$ be a volume of $G$.
 In this paper we prove that
 $$
      \limsup_{|G| \to \infty}   (\det \Gamma \#(\Gamma\cap G)-|G|)/\ln^{s-1} |G| >0.
 $$
 Thus the known estimate $\det \Gamma  \#(\Gamma\cap G)=|G| +O(\ln^{s-1} |G|)$    is exact.
  We obtain also a similar result for the low discrepancy sequence corresponding to $\Gamma$.
\end{abstract}
%
%
\section{Introduction.}


{\bf 1.1. Lattice points.}
  Let $ \Gamma \subset \RR^s $ be a lattice, i.e., a
 discrete subgroup of $\RR^s $
 with a compact fundamental set $ \RR^s/\Gamma$, 
$\det \Gamma$=vol$(\RR^s/\Gamma)$.  
 Let $N_1,...,N_s>0$ be reals, $\bN = (N_1,...,N_s)$,  $B_{\bN} =[0,N_1)\times \cdots \times [0,N_s) $, 
 $\vol( B_{\bN} )$ the volume of $B_{\bN}$, $tB_{\bN} $ the dilatation
  of $ B_{\bN} $ by a factor $ t>0$,  $tB_{\bN}+\bx$  the translation of 
 $tB_{\bN} $ by  a vector $ \bx \in \RR^s$,  $(x_1,...,x_s) \cdot (y_1,...,y_s) =(x_1y_1,...,x_s y_s)$, and let $(x_1,...,x_s) \cdot B_{\bN} =\{ (x_1,...,x_s) \cdot (y_1,...,y_s) \; | \;
 (y_1,...,y_s) \in B_{\bN} \} $. 
Let
\begin{equation}\label{1.0}
  \cN (B_{\bN} + \bx ,\Gamma)=\#(B_{\bN} + \bx \cap \Gamma )=\sum_{\bgamma\in
  \Gamma} \d1_{B_{\bN} +\bx} (\bgamma)
\end{equation}
be the number of points of the lattice $ \Gamma $ lying inside the
parallelepiped $B_{\bN} $, where we denote by $ \d1_{B_{\bN} +\bx} (\bgamma)$ the indicator function of $ B_{\bN} +\bx$. 
We define the error $\cR(B_{\bN}+\bx,\Gamma)$  by setting
\begin{equation}\label{1.2}
 \cN (B_{\bN}+\bx,\Gamma)= (\det \Gamma)^{-1} \vol( B_{\bN})   \;+\;\cR(B_{\bN}+\bx,\Gamma).
\end{equation}
Let  ${\rm Nm} ( \bx)= x_1 x_2 \ldots x_s$  for $\bx=(x_1, \dots, x_s) $.
 The lattice $ \Gamma \subset \RR^s $ is {\it admissible} if
\begin{equation} \nonumber
  {\rm Nm} \;\Gamma =\inf_{\bgamma \in \Gamma \setminus \{0\}} | {\rm Nm} (\bgamma)| >0.
\end{equation}
Let $ \Gamma $  be an admissible lattice. In 1994, Skriganov [Skr] proved the following theorem:
\\

{\bf Theorem A.} {\it Let $\bt =(t_1,...,t_s)$. Then 
\begin{equation}\label{1.6}
   | \cR(\bt \cdot  [-1/2,1/2)^s +\bx, \Gamma) |  <c_0 (\Gamma) \log_2^{s-1} (2+ |{\rm Nm} (\bt)|),
\end{equation}
where the constant  $ c_0(\Gamma)$  depends upon the lattice
$ \Gamma $ only by means of the invariants $ \det \Gamma $ and ${\rm Nm} \; \Gamma $.}\\

 In [Skr, p.205], Skriganov conjectured that the bound (\ref{1.6})  is the best possible. 
 In this paper we prove this conjecture. 

Let $\cK$ be a totally real algebraic number field of degree $s \ge 2$, and let 
 $\sigma$ be the canonical embedding of $\cK$ in the Euclidean space
 $\RR^s$, $ \sigma : \cK \ni \xi \rightarrow \sigma (\xi) =
 (\sigma_1(\xi), \ldots, \sigma_s (\xi)) \in \RR^s $, where $
 \{\sigma_j \}_{j=1}^s $ are $s$ distinct embeddings of $\cK$ in the field $
 \RR$ of real numbers. Let  $N_{\cK/\QQ}(\xi)$ be the norm of  $\xi \in \cK$. By [BS, p. 404],
\begin{equation} \nonumber
    N_{\cK/\QQ}(\xi)  = \sigma_1 (\xi) \cdots \sigma_s (\xi), \quad {\rm and}  \quad    |N_{\cK/\QQ}(\alpha)| \geq 1
\end{equation}
for all algebraic integers  $\alpha \in \cK \setminus \{ 0 \}$. We see that
$ |{\rm Nm} (\sigma (\xi))|=  | N_{\cK/\QQ}(\xi)|$. Let $\cM$  be a  full $\ZZ$ module  in $ \cK$ and let $\Gamma_{\cM}$ be the lattice 
 corresponding to $\cM$ under the embedding $\sigma$. 
Let $(c_{\cM})^{-1} >0$ be an integer such that $(c_{\cM})^{-1} \gamma$ are algebraic integers for all $ \gamma \in \cM $. Hence
\begin{equation}\nonumber
   {\rm Nm} \;\Gamma_{\cM}   \geq c_{\cM}^s.
\end{equation} 
Therefore, $\Gamma_{\cM} $   is an admissible lattice.
In the following, we will use notations $\Gamma=\Gamma_{\cM}$, and $N =N_1N_2...N_s \geq 2$. 
 In $\mathsection 2$ we will prove  the following theorem:\\

{\bf Theorem 1.} {\it With the  above notations, there exist $c_1(\cM)>0$ such that
\begin{equation}\label{1.12}
      \sup_{\btheta \in [0,1]^s}   |\cR( B_{\btheta  \cdot \bN}  +\bx, \Gamma_{\cM})  |\geq c_1({\cM})\log_2^{s-1} N
\end{equation}  
  for all $\bx \in \RR^{s}$.}\\

In [La1, Ch. 5], Lang considered the  lattice point problem in the adelic setting. In [La1] and [NiSkr], the upper bound for the  lattice point remainder problem in parallelotopes was found.
 In a forthcoming paper, we will prove that the lower bound (\ref{1.12}) can be extended to the adelic case. Namely, we will prove that the upper bound in [NiSkr] is exact for the case of totally real algebraic number fields.
\\

{\bf 1.2. Low discrepancy sequences.}
 Let $((\beta_{k,N})_{k=0}^{N-1})$ be a $N$-point set in an $s$-dimensional unit cube $[0,1)^s$,
$B_{\by}=[0,y_1) \times \cdots \times [0,y_s) \subseteq [0,1)^s $,
\begin{equation}\label{1.14}
\Delta(B_{\by}, (\beta_{k,N})_{k=0}^{N-1}  )= \#\{0 \leq k
<N  \;|\; \beta_{k,N}\in B_{\by}\}-Ny_1 \ldots y_s.
\end{equation}
We define the star {\it discrepancy} of a 
$N$-point set $(\beta_{k,N})_{k=0}^{N-1}$ as
\begin{equation} \label{1.16}
   \emph{D}^{*}(N)=\emph{D}^{*}((\beta_{k,N})_{k=0}^{N-1}) = 
    \sup_{ 0<y_1, \ldots , y_s \leq 1} \; | \frac{1}{  N}
  \Delta(B_{\by},(\beta_{k,N})_{k=0}^{N-1}) |.
\end{equation}

In 1954, Roth proved that there exists a constant $ \dot{c}_1>0 $, such
that
\begin{equation} \nonumber
N\emph{D}^{*}((\beta_{k,N})_{k=0}^{N-1})>\dot{c}_1(\ln N)^{s-1\over 2},   
\end{equation}
for all $N$-point sets $(\beta_{k,N})_{k=0}^{N-1}$.

 \texttt{ Definition 1.} {\it  A sequence  of point sets $((\beta_{k,N})_{k=0}^{N-1})_{N=1}^{\infty}$ is of 
 \texttt{low discrepancy} (abbreviated
l.d.p.s.) if $ \emph{D}^{*}((\beta_{k,N})_{k=0}^{N-1})=O(N^{-1}(\ln
N)^{s-1}) $ for $ N \rightarrow \infty $. }

 For examples of  l.d.p.s. see e.g. in [BC], [DrTi], and [Skr].
Consider a lower bound for    l.d.p.s. 
According to the well-known conjecture (see, e.g., [BC, p.283]), there exists 
a constant $ \dot{c}_2>0 $, such that
\begin{equation}   \label{1.22}
  {N\emph{D}^{*}((\beta_{k,N})_{k=0}^{N-1}) >\dot{c}_2 (\ln N)^{s-1}} 
\end{equation}
for all $N$-point sets $(\beta_{k,N})_{k=0}^{N-1}$.   
In 1972, W. Schmidt proved this conjecture for $ s=2 $.
In 1989,   Beck  [Be] proved that  $N\emph{D}^{*}(N) \geq \dot{c} \ln N  (\ln\ln N)^{1/8-\epsilon}$  for $s=3$ and some $\dot{c}>0$. In 2008, Bilyk,  Lacey
and Vagharshakyan (see [Bi, p.147], [BiLa, p.2]), proved   in all
 dimensions $s \geq 3$ that there exists some $\dot{c}(s), \eta >0$ for which the following
estimate holds for all  $N$-point sets :
$N\emph{D}^{*} (N)>\dot{c}(s)(\ln N)^{\frac{s-1}{2}  +\eta}$.

There exists another conjecture  on the lower bound for the discrepancy function:
 there exists 
a constant $\dot{c}_3>0 $, such that
\begin{equation}   \nonumber
  {N\emph{D}^{*}((\beta_{k,N})_{k=0}^{N-1}) >\dot{c}_3 (\ln N)^{s/2}} 
\end{equation}
for all $N$-point sets $(\beta_{k,N})_{k=0}^{N-1}$  (see [Bi, p.147], [BiLa, p.3] and [ChTr, p.153]). 

 Let $\cW= (\Gamma +\bx) \cap [0,1)^{s-1} \times [0,\infty)$. We enumerate $\cW$ by the sequence $(z_{1,k}(\bx), z_{2,k}(\bx))$  with
  $z_{1,k}(\bx) \in [0,1)^{s-1} $ and  $z_{2,k}(\bx) \in [0,\infty) $. In [Skr], Skriganov proved that the point set 
   $((\beta_{k,N}(\bx))_{k=0}^{N-1})$ with
    $\beta_{k,N}(\bx) =(z_{1,k}(\bx), z_{2,k}(\bx)/N) $ is of low discrepancy (see also [Le1]).  In $\mathsection 2.10$  we will prove \\

{\bf Theorem 2.} {\it With the notations as above, there exist $c_2(\cM)$ such that
\begin{equation}\label{1.26}
      \sup_{y_s\in [0,1]}  \emph{D}^{*}((\beta_{k,N}(\bx))_{k=0}^{[y_s N]})   \geq c_2(\cM)\log_2^{s-1} N
\end{equation}  
  for all $\bx \in \RR^{s}$.}

This result  support the conjecture (\ref{1.22}).
In a forthcoming paper we will prove that  (\ref{1.26}) is also true for Halton's sequence, $(t,s)$-sequence, and l.d.p.s. from [Le2].
\\

\section{Proof of Theorems.}
{\bf 2.1. Poisson summation formula.} 
  It is known that the set $\cM^{\bot}$ of all
  $\beta \in \cK$, for which $ {\rm Tr}_{\cK/  \QQ}(\alpha \beta)  \in\ZZ$ for all $\alpha \in \cM$, is also a full $\ZZ$
 module ({\it the dual of the module $\cM$}) of the field $K$ (see [BS], p. 94).
 Recall that  the dual lattice  $ \Gamma_\cM^\bot $ consists of all vectors  $
\bgamma^\bot \in \RR^{s}$ such that the inner product $<\bgamma^\bot,\bgamma>$
belongs to $\ZZ $ for each $\bgamma \in \Gamma $.  Hence $\Gamma_{\cM^{\bot}} =\Gamma_\cM^{\bot}$.
  Let $\fO$ be the ring of integers of the field $\cK$, and let $a \cM^{\bot} \subseteq \fO  $ for some $a \in \ZZ \setminus 0$. 
  By (\ref{1.0}),   we have $ \cN(B_{\bN}  +\bx  ,\Gamma_{\cM}  ) = \cN( a^{-1}B_{\bN}  +a^{-1}\bx  ,\Gamma_{a^{-1}\cM}  )$. Therefore, to prove Theorem 1 it suffices
consider only the case $\cM^{\bot} \subseteq \fO $. We set
\begin{equation}\label{2.0}
   p_1 = \min\{ b \in \ZZ \; | \; b \fO \subseteq \cM^{\bot} \subseteq \fO , \;\; b>0  \}.
\end{equation}  
We will use the same notations for elements of $\fO$ and $\Gamma_{\fO}$. 
  Let $\cD_{\cM}$ be the ring of coefficients of the full module $\cM$,  $\cU_{\cM}$ be the group of units of $\cD_M$,
and let   $\eta_{1},..., \eta_{s-1}$ be the set of fundamental units of $\cU_{\cM}$. According to the Dirichlet
 theorem (see e.g., [BS, p. 112]), every unit $\epsilon \in \cU_{\cM}$ has a unique representation in the form
\begin{equation}\label{2.2}
   \epsilon = (-1)^a\eta_{1}^{a_1} \cdots \eta_{s-1}^{a_{s-1}},  
\end{equation} 
where $a_1,...,a_{s-1}$ are rational integers and $a \in \{0,1\}$.   
 It is easy to proof (see e.g. [Le3, Lemma 1]) that there exists a constant $c_3>1$ such that  for all $\bN$ there exists  $\eta(\bN) \in \fU_{\cM}$ with 
$  |N_i^{'} N^{-1/s}| \in [1/c_3, c_3]$, 
where $ N_i^{'} =  N_i |\sigma_i(\eta(\bN))|$, $ i=1,...,s $, and $N=N_1 \cdots N_s$.  
Let $\sigma(\eta(\bN)) = ( \sigma_1(\eta(\bN)),...,\sigma_s(\eta(\bN)))$.
We see that $\sigma(\eta(\bN)) \cdot (\btheta \cdot B_{\bN} +\bx) = \btheta \cdot B_{\bN^{'}} +\bx_1$ and
\begin{equation}\nonumber
         \bgamma \in \Gamma_\cM     \cap  (\btheta \cdot B_{\bN} +\bx)  \Leftrightarrow \bgamma \cdot \sigma(\eta(\bN))  \in \Gamma_\cM   
           \cap  (\btheta \cdot B_{\bN^{'}} +\bx_1 )),
\end{equation}  
with $\bx_1= \sigma(\eta(\bN) \cdot \bx  +  \sigma(\eta(\bN)) \cdot \bN/2 -\bN^{'}/2$. 
Hence
\begin{equation} \nonumber
   \cN(\btheta \cdot B_{\bN} +\bx  ,\Gamma_\cM  )  = \cN(\btheta \cdot B_{\bN^{'}} +\bx_1,\Gamma_\cM  ).
\end{equation}  
By (\ref{1.2}), we have 
\begin{equation}\nonumber
   \cR(\btheta \cdot B_{\bN} +\bx  ,\Gamma_\cM  )  = \cR(\btheta \cdot B_{\bN^{'}} +\bx_1 ,\Gamma_\cM  ).
\end{equation}  
Therefore, without loss of generality, we can assume that
\begin{equation}\label{2.8}
    N_i N^{-1/s} \in [1/c_3, c_3], \quad  i=1,...,s.  
\end{equation}  
Note that in this paper $O$-constants and constants $c_1,c_2,...$ depend only on $\cM$.

We shall need the Poisson summation formula: 
\begin{equation}\label{2.32}
             \det \Gamma   \sum_{\bgamma\in \Gamma} f(\bgamma-X)  =   \sum_{\bgamma\in \Gamma^\bot }  \widehat{f}(\bgamma)  e(\langle\bgamma,\bx\rangle),
\end{equation}
where  
\begin{equation} \nonumber
            \widehat{f}(Y) = \int_{\RR^s}{ f(X) e(\langle \by,\bx\rangle) d\bx}  
\end{equation}
is the Fourier transform of $f(X)$, and $e(x) =exp(2\pi \sqrt{-1}x), \langle\by,\bx\rangle = y_1x_1+ \cdots +y_s x_s$. Formula (\ref{2.32}) holds for functions $f(\bx)$ with period lattice $\Gamma$ 
if one of the functions $f$ or $\widehat{f}$ is integrable  and belongs to the class $C^{\infty}$ (see e.g.  [StWe, p. 251]).

Let $ \widehat{\d1}_{B_{ \bN}}(\bgamma)$ be the Fourier transform of the indicator function $\d1_{B_{ \bN}}(\bgamma)$.
It is easy to prove that  $\widehat{\d1}_{B_{ \bN}}(\bs) = N_1\cdots N_s$ and
\begin{equation}\label{2.34b}
 \widehat{\d1}_{B_{ \bN}}(\bgamma) = \prod_{i=1}^s \frac{e(  N_i \gamma_i) -1}{2\pi \sqrt{-1}\gamma_i}
    = \prod_{i=1}^s \frac{\sin(\pi  N_i  \gamma_i)}{ \pi \gamma_i} e\Big(\sum_{i=1}^s   N_i  \gamma_i/2   \Big) \; {\rm for}\; \Nm (\bgamma \neq 0).
\end{equation}

We fix a nonnegative even function $ \omega(x),\; x \in \RR,$ of the
class $ C^\infty $, with a support  inside the segment $ [-1/2,  1/2] $, and satisfying the condition $ \int_{\RR }\omega(x)dx=1$.
We set $\Omega(\bx) =  \omega( x_1)   \cdots   \omega( x_s)  $,    $\Omega_{\tau}(\bx) = \tau^{-s} \Omega(\tau^{-1} x_1,...,\tau^{-1} x_s)  $, $\tau>0$, and
\begin{equation}\label{2.10}
 \widehat{\Omega} (\by)   =        \int_{\RR^s }e(\langle \by,\bx\rangle) \Omega(\bx)d\bx.
\end{equation}
 Notice that the Fourier transform $\widehat{\Omega}_{\tau}( \by) =\widehat{\Omega} (\tau \by)$ of the
function $ \Omega_{\tau}(\by)$ satisfies the bound
\begin{equation}\label{2.12}
    |\widehat{\Omega} (\tau \by)|<\dot{c}(s, \omega)(1+ \tau|\by|)^{-2s} .
\end{equation}
It is easy to see that
\begin{equation}\label{2.13}
    \widehat{\Omega}(\by)   =\widehat{\Omega}(\bs) +O( |\by|) =1+O( |\by|) \quad {\rm for} \quad   |\by| \to 0.
\end{equation} \\

 {\bf Lemma 1.} {\it There exists a constant $c>0$, such that we have for $  N> c$ 
\begin{equation}\nonumber
   |\cR (B_{\btheta \cdot\bN}+\bx,\Gamma ) -  \ddot{\cR} (B_{\btheta \cdot\bN}+\bx,\Gamma )| 
    \leq     2^s,
\end{equation} 
where}
\begin{equation}\label{2.34a}
    \ddot{\cR} (B_{\btheta \cdot\bN}+\bx,\Gamma) = (\det \Gamma)^{-1}  \sum_{\bgamma\in \Gamma^\bot \setminus  \{0\}}
   \widehat{\d1}_{B_{\btheta \cdot\bN}}(\bgamma)\widehat{\Omega}
    (\tau  \bgamma)e(\langle\bgamma,\bx\rangle), \quad \tau =N^{-2}.
\end{equation}
{\bf Proof.} 
Let $B^{\pm \tau}_{\btheta \cdot \bN} =  [0, \max(0,\theta_1N_1\pm \tau)) \times \cdots    \times   [0, \max(0,\theta_sN_s\pm \tau))$, and let $\d1_B(x)$ be the indicator function of $B$. 
 We consider the convolutions of the functions $\d1_{B^{\pm \tau}_{\btheta \cdot \bN}}( \bgamma)$  and $\Omega_{\tau}(\by)$ : 
\begin{equation}\label{2.14}
    \Omega_{\tau} \ast \d1_{B^{\pm \tau}_{\btheta \cdot \bN}}  (\bx) =   \int_{\RR^s}   \Omega_{\tau} (\bx-\by)  \d1_{B^{\pm \tau}_{\btheta \cdot \bN}} (\by) d\by . 
\end{equation}
It is obvious that the nonnegative functions (\ref{2.14}) are of class $ C^\infty $ and are compactly supported in 
$\tau$-neighborhoods of the bodies $B^{\pm \tau}_{\btheta \cdot \bN} $, respectively. 
We obtain
\begin{equation}\label{2.16}
     \d1_{B_{\btheta \cdot \bN}^{-\tau}} (\bx)   \leq \d1_{B_{\btheta \cdot \bN}} (\bx)  \leq  \d1_{B_{\btheta \cdot \bN}^{+\tau}} (\bx) , \;\;
     \d1_{B_{\btheta \cdot \bN}^{-\tau}} (\bx)  \leq \Omega_{\tau}  \ast  \d1_{B_{\btheta \cdot \bN} }(\bx)     \leq   \d1_{B_{\btheta \cdot \bN}^{+\tau}} (\bx) .
\end{equation}
Replacing $\bx$ by $\bgamma -\bx$ in (\ref{2.16}) and summing these inequalities over $\bgamma \in\Gamma = \Gamma_\cM$, we find from (\ref{1.0}), that
\begin{equation}\nonumber
  \cN (B^{-\tau}_{\btheta \cdot\bN}+\bx,\Gamma ) \leq    \cN (B_{\btheta \cdot \bN}+\bx,\Gamma ) \leq   \cN (B^{+\tau}_{\btheta \cdot\bN}+\bx,\Gamma ),
\end{equation}
and
\begin{equation}\nonumber
         \cN (B^{-\tau}_{\btheta \cdot\bN}+\bx,\Gamma )    \leq  \dot{\cN} (B_{\btheta \cdot\bN}+\bx,\Gamma )   \leq
                        \cN (B^{+\tau}_{\btheta \cdot\bN}+\bx,\Gamma ),     
\end{equation}
where
\begin{equation}\label{2.20}
         \dot{\cN} (B_{\btheta \cdot\bN}+\bx,\Gamma ) = \sum_{\bgamma\in \Gamma}  \Omega_{\tau}  \ast \d1_{B_{\btheta \cdot \bN}} (\bgamma - \bx).
\end{equation}
Hence
\begin{equation}\nonumber
   -  \cN (B^{+ \tau}_{\btheta \cdot\bN}+\bx,\Gamma ) + \cN (B^{- \tau}_{\btheta \cdot\bN}+\bx,\Gamma ) 
\end{equation}
\begin{equation}\nonumber
     \leq  \dot{\cN} (B_{\btheta \cdot\bN}+\bx,\Gamma ) - \cN (B_{\btheta \cdot\bN}+\bx,\Gamma )  \leq  \cN (B^{+ \tau}_{\btheta \cdot\bN}+\bx,\Gamma ) - \cN (B^{- \tau}_{\btheta \cdot\bN}+\bx,\Gamma ) .
\end{equation}
Thus
\begin{equation}\label{2.22}
       |\cN (B_{\btheta \cdot\bN}+\bx,\Gamma ) -  \dot{\cN} (B_{\btheta \cdot\bN}+\bx,\Gamma )| \leq  \cN (B^{+ \tau}_{\btheta \cdot\bN}+\bx,\Gamma ) - \cN (B^{- \tau}_{\btheta \cdot\bN}+\bx,\Gamma ) .
\end{equation}
Consider the right side of this inequality.
We have that $B^{+ \tau}_{\btheta \cdot\bN}  \setminus B^{- \tau}_{\btheta \cdot\bN}$ is the union of boxes $B^{(i)}, \; i=1,...,2^s-1$, where
\begin{equation}  \nonumber
 \vol(B^{(i)}) \leq  \vol( B_{\bN}^{+\tau})  -  \vol(B_{\bN}^{-\tau})   \leq \prod_{i=1}^s (N_i + \tau) - \prod_{i=1}^s (N_i - \tau) 
\end{equation}
\begin{equation}  \nonumber
 \leq N \Big( \prod_{i=1}^s (1 + \tau) - \prod_{i=1}^s (1 - \tau) \Big) <\ddot{c}_s N \tau =\ddot{c}_s  / N,  , \quad \tau =N^{-2},
\end{equation} 
with some $\ddot{c}_s>0$.  From (\ref{2.0}), we get $\cM \supseteq p_1^{-1}\fO$. Hence $|{\rm Nm}( \bgamma)| \geq p_1^{-s}$ for $\bgamma \in \Gamma_{\cM} \setminus \bs$.
We see that $ |{\rm Nm}( \bgamma_1 - \bgamma_2  )| \leq  \vol(B^{(i)} +\bx) < p_1^{-s} $  for $ \bgamma_1, \bgamma_2  \in B^{(i)} +\bx $
 and  $  N>\ddot{c}_s  p_1^{s} $.  Therefore, the
  box $B^{(i)}+\bx $ 
contains at most one point of $\Gamma_\cM$ for $  N>\ddot{c}  p_1^{s}  $. By (\ref{2.22}), we have
\begin{equation}\label{2.28}
    |\dot{\cN} (B_{\btheta \cdot\bN}+\bx,\Gamma ) - \cN (B_{\btheta \cdot\bN}+\bx,\Gamma ) |  \leq  2^s -1,  \quad {\rm for}  \quad  N>\ddot{c}  p_1^{s} .
\end{equation}
Let
\begin{equation}\label{2.30}
   \dot{\cR} (B_{\btheta \cdot\bN}+\bx,\Gamma ) = \dot{\cN} (B_{\btheta \cdot\bN}+\bx,\Gamma ) - \frac{ \vol(B_{\btheta \cdot\bN})}{\det \Gamma}   .
\end{equation}  

By  (\ref{2.20}), we obtain that $ \dot{\cN} (B_{\btheta \cdot\bN}+\bx,\Gamma )$ is a periodic
function of $ \bx \in \RR^n $ with the period lattice $
\Gamma$. 
Applying the Poisson summation formula to the series (\ref{2.20}), and bearing in mind that $ \widehat{\Omega}_{\tau} (\by) = \widehat{\Omega} (\tau \by)$,  
 we get from (\ref{2.34a})
\begin{equation}  \nonumber
\dot{\cR} (B_{\btheta \cdot\bN}+\bx,\Gamma) = \ddot{\cR} (B_{\btheta \cdot\bN}+\bx,\Gamma) .
\end{equation}
Note that (\ref{2.12})
ensure the absolute convergence of the series (\ref{2.34a}) over $\bgamma\in \Gamma^\bot \setminus  \{0\}$.
Using (\ref{1.2}), (\ref{2.28}) and (\ref{2.30}) , we obtain the assertion of Lemma 1.   \qed \\

Let $\eta(t)=\eta(|t|)$, $t \in \RR^1$ be an even function of the class $C^{\infty}$; moreover, let $\eta(t)=0$ for $|t| \leq 1$, $0 \leq \eta(t) \leq 1$ for $|t| \leq 2$  and $\eta(t)=1$ for $|t| \geq 2$.
Let $n=s^{-1}\log_2 N $, $M=[\sqrt{n}]$ , and
\begin{equation}\label{2.35a}
              \eta_M(\bgamma) = 1-\eta(2|\Nm(\bgamma)|/M).  
\end{equation}
By  (\ref{2.34b}) and (\ref{2.34a}), we have
\begin{equation}\label{2.35}
 \dot{\cR} (B_{\btheta \cdot\bN}+\bx,\Gamma )= (\pi^s  \det \Gamma)^{-1} (\cA(\bx,M) +\cB(\bx,M)),    
\end{equation}
where
\begin{equation}\nonumber
 \cA(\bx,M) = \sum_{\bgamma\in \Gamma^\bot \setminus {\bs} } 
          \prod_{i=1}^s \sin(\pi \theta_i N_i \gamma_i) \frac{    \eta_M(\bgamma) \widehat{\Omega} (\tau  \bgamma)e(\langle\bgamma,\bx\rangle +\dot{x})}{\Nm(\bgamma)},
\end{equation}
\begin{equation}\nonumber
    \cB(\bx,M) = \sum_{\bgamma\in \Gamma^\bot \setminus {\bs} } 
          \prod_{i=1}^s \sin(\pi \theta_i N_i \gamma_i) \frac{   (1- \eta_M(\bgamma)) \widehat{\Omega} (\tau  \bgamma)e(\langle\bgamma,\bx\rangle +\dot{x})}{\Nm(\bgamma)},
\end{equation}
 with $\dot{x} = \sum_{1 \leq i  \leq s}  \theta_i N_i  \gamma_i/2  $. Let
\begin{equation}\nonumber
        \bE(f) =\int_{[0,1]^s} f(\btheta)d\btheta.      
\end{equation} 
By the triangle inequality, we get
\begin{equation}\label{2.46}
   \pi^s  \det \Gamma     \sup_{\btheta \in [0,1]^s} |  \dot{\cR} (B_{\btheta \cdot\bN}+\bx,\Gamma )  |  \geq |\bE(\cA(\bx,M))|     - |\bE(\cB(\bx,M))|.
\end{equation} 
   In $\mathsection 2.5$ we will find the lower bound of $|\bE(\cA(\bx,M))|$ and   in $\mathsection 2.9$ we will find the upper bound of $|\bE(\cB(\bx,M))|$. 
\\
\\

{\bf 2.2. The  logarithmic space and the fundamental domain.} We consider Dirichlet's Unit Theorem (\ref{2.2}) applied to the ring of integers $\fO$. Let
 $\bepsilon_1,..., \bepsilon_{s-1}$ be the set of fundamental units of $\cU_{\fO}$.
We set $l_i(\bx) = \ln|x_i|$, $i=1,...,s$, $\bl(\bx) =(l_1(\bx),...,l_s(\bx))$, $\b1=(1,...,1)$, where
  $\bx \in \RR^s$  and $\Nm(\bx) \neq \bs $.
By [BS, p. 311],  the set of vectors $\b1, \bl(\epsilon_1),..., \bl(\epsilon_s-1))$  is a basis for $\RR^s $. Any vector $\bl(\bx) \in \RR^s $ ($\bx \in \RR^s, \; \Nm(\bx) \neq \bs $) can be 
represented in the form
\begin{equation}\label{3.0}
       \bl(\bx) = \xi \b1 + \xi_1 \bl(\bepsilon_1) + \cdots +  \xi_{s-1} \bl(\bepsilon_{s-1})  
\end{equation}
where $\xi,\xi_1,...,\xi_{s-1}$ are real numbers.  In the following we will need the next definition

    \texttt{ Definition 2.} [BS, p. 312] {\it 
 A subset $\cF$ of the space $\RR^s$ is called a  \texttt{fundamental domain}
for the field $\cK$ if it consists of all points $\bx$ which satisfy the following 
conditions: 
 $\Nm(\bx) \neq \bs$, in the representation (\ref{3.0}) the coefficients $\xi_i$ $(i = 1, ... , s-~1~)$ satisfy the
inequality $0 \leq \xi_i < 1$, $x_1 >0$.} \\

{\bf Theorem B}. [BS, p. 312] In every class of associate numbers $( \neq 0)$ of the field $\cK$, there
is one and only one number whose geometric representation in the space $\RR^s$
lies in the fundamental domain $\cF$ . \\

{\bf Lemma A.} [Wi, p.59, Theorem 2, ref. 3] {\it Let $\dot{\Gamma} \subset \RR^k$ be a lattice, $\det \dot{\Gamma} =1$, $\cQ \subset \RR^k$ a compact convex body
 and $r$ the radius of
its greatest sphere in the interior. Then
\begin{equation}\nonumber
       \vol( \cQ) \Big(1- \frac{\sqrt{k}}{2r}\Big)  \leq  \# \dot{\Gamma}  \cap \cQ  \leq  \vol( \cQ )\Big(1+ \frac{\sqrt{k}}{2r}\Big) ,
\end{equation}  
provided $r> \sqrt{k}/2$.} \\

Let $\dot{\Gamma} \subset \RR^k$ be an arbitrary lattice. We derive from Lemma A
\begin{equation}\label{3.2}
      \sup_{\bx \in \RR^s} | \#\dot{\Gamma} \cap (t\cQ +\bx) -  t^k \vol( \cQ) /\det\dot{\Gamma} | = O(t^{k-1}) \quad {\rm for} \quad t \to \infty.
\end{equation}  
See also [GrLe, p. 141,142].\\

{\bf Lemma 2.} {\it  Let $ \bepsilon^{\bk}_{max} =\max_{1 \leq i \leq s}|(\bepsilon^{\bk})_i|$ and 
  $\bepsilon^{\bk}_{min} = \min_{1 \leq i \leq s}|(\bepsilon^{\bk})_i|$. There exists a constant $c_4,c_5>0$,  such that
\begin{equation}\label{3.4}
       \# \{ \bk \in \ZZ^{s-1} \; | \; \bepsilon^{\bk}_{max} \leq  e^{t} \} =  c_4t^{s-1}  + O(t ^{s-2} )
\end{equation} 
and
\begin{equation}\label{3.4a}
       \# \{ \bk \in \ZZ^{s-1} \; | \; \bepsilon^{\bk}_{min} \geq  e^{-t}  \} =  c_5t^{s-1}  + O(t ^{s-2} ) .
\end{equation} }
{\bf Proof}.  By (\ref{3.0}),  we have that the left hand sides of (\ref{3.4}) and (\ref{3.4a}) are equal to 
\begin{equation}\nonumber 
       \# \{ \bk \in \ZZ^{s-1} \; | \;  \sum_{i=1}^{s-1} k_i l_j(\bepsilon_i) \leq t, \; j=1,...,s  \}, 
\end{equation}  
and 
\begin{equation}\nonumber 
       \# \Big\{ \bk \in \ZZ^{s-1} \; | \;  \sum_{i=1}^{s-1} k_i l_j(\bepsilon_i) \geq -t, \; j=1,...,s  \Big\}, 
\end{equation}  
respectively. Let
\begin{equation}\nonumber 
   \cQ_1       =  \{ \bx \in \RR^{s-1} |  \dot{x_j}  \leq 1,\;   j \in [1,s]  \} \; \; {\rm and} \;\; \cQ_2       =  \{ \bx \in \RR^{s-1}  |  \dot{x_j}  \geq -1,\; j \in [1,s]\},
\end{equation}
with $\dot{x_j} = x_1 l_j(\bepsilon_1) + \cdots + x_{s-1} l_j(\bepsilon_{s-1})  $.   
We see $\dot{x_1}+\cdots+\dot{x_s} =0$. Hence $ \dot{x_j}  \geq -s+1$ for $\bx \in \cQ_1 $ and $ \dot{x_j}  \leq s-1$ for $\bx \in \cQ_2 $ $(j=1,...,s )$. 
  By [BS, p. 115], we get $ \det (l_i(|\bepsilon_j|)_{i,j=1,...,s-1})  \neq 0$.
 Hence, $  \cQ_i  $ is the compact convex set in $\RR^{s-1}$, $i=1,2$.
Applying (\ref{3.2}) with $k=s-1$, and $\dot{\Gamma}  =\ZZ^{s-1}$, we obtain  the assertion of  Lemma~2. \\ 
\hfill \qed  \\

Let $ cl(\cK)$ be the ideal class group of $\cK$, $h = \# cl(\cK)$, and $ cl(\cK) = \{C_1,..., C_h \}$. In the ideal class $C_i$, we choose an integral 
 ideal $\fa_i$, $i=1,...,h$.  Let $\fN(\fa)$ be the absolute norm of ideal $\fa$. If $h=1$, then we set $p_2=1$ and $\Gamma_1= \Gamma_{\fO}$. 
 Let $h>1$, $i \in [1,h]$,
\begin{equation}\label{3.5}
      \cM_i =\{ u \in \fO \; | \; u \equiv 0 \mod \fa_i   \}, \quad \Gamma_i = \sigma(\cM_i),  \;\; {\rm and}\;\;
       p_2 = \prod_{i=1}^h  \fN (\fa_i) .  
\end{equation} \\

{\bf Lemma 3.} {\it Let $w \geq 1$, $i \in [1,h]$, $\FF_{M_1}(\bvarsigma) =\{\by \in \cF \;|\; |\Nm (\by)| <M_1, \\ \sgn(y_i) =\varsigma_i,  i=1,...,s   \}$, where   $\sgn(y) = y/|y|$ for $y \neq 0$
 and $\bvarsigma=(\varsigma_1,...,\varsigma_s) \in \{ -1, 1\}^s$. Then there exists $c_{6,i}>0$,  such that}
\begin{equation}\nonumber
  \sup_{\bx \in \RR^s}    \Big|  \sum_{ \gamma \in (w \Gamma_{i} + \bx )\cap \FF_{M_1}(\bvarsigma)} 1  - c_{6,i} M_1/w^s \Big| =  O(M_1^{1-1/s}), \quad {\rm for} \quad M_1 \to \infty.
\end{equation} 
{\bf Proof}. It is easy to see that $\FF_{M_1}(\bvarsigma) ={M_1}^{1/s}\FF_1(\bvarsigma)$. By  [BS, p. 312 ], the fundamental domain $\cF$ is  a cone in $\RR^s$.
  Let  $\dot{\FF} =\{ \by \in \cF |\;|y_i| \leq y_0, \sgn(y_i) =\varsigma_i,\\ i=1,...,s \}$
 and let  $\ddot{\FF} =\{ \by \in \dot{\FF}\;|\; |\Nm(\by)| \geq 1\}$, where  $y_0 =\sup_{\by \in \FF_1(\bvarsigma), i=1,...,s} |y_i|$.
 We see  that $\FF_1(\bvarsigma) = \dot{\FF} \setminus \ddot{\FF} $ and $\dot{\FF}, \ddot{\FF} $ are compact convex sets.
 Using (\ref{3.2}) with $k=s$, $\dot{\Gamma} =w \Gamma_{i} $,  and $t={M_1}^{1/s}$, we obtain  the assertion of  Lemma 3. \qed  \\
\\

{\bf 2.3. Construction a Hecke character by using Chevalley's theorem.} 
Let $\fm$ be an integral ideal of the number field   $\cK$, and let $\cJ^{\fm}$ be the group
of all ideals of $\cK$ which are relatively prime to $\fm$. Let $S^1=\{z \in \CC \; | \; |z| = 1\} $.

    \texttt{ Definition 3.} [Ne, p. 470]  { \it A \texttt{ Hecke character} $\mod \fm$ is a character $\chi :\; \cJ^{\fm} \to S^1$
for which there exists a pair of characters
\begin{equation}\nonumber
          \chi_f :\;   (\fO/\fm)^{*}   \to S^1,    \qquad       \chi_{\infty} :\; (\RR^{*})^s \to S^1,  
\end{equation}
such that
\begin{equation}\nonumber
       \chi((a)) = \chi_f(a)   \chi_{\infty} (a)   
\end{equation}
for every algebraic integer $a \in \fO$ relatively prime to $\fm$.  }


The character taking the value one for all group elements will be called the trivial character.\\

    \texttt{Definition 4.} {\it 
Let $A_1,...,A_d$ be  invertible $s\times s$ commuting  matrices with integer entries. 
A sequence of  matrices $A_1,...,A_d$ 
is said to be { \texttt{ partially hyperbolic}} if for all $(n_1,...,n_d) \in \ZZ^d \setminus \{0 \}$ 
 none of the eigenvalues of $A_1^{n_1}...A_d^{n_d}$ are roots of unity.
}\\

 We need the following variant of Chevalley's theorem ([Ch], see also [Ve]):

{\bf Theorem C} [KaNi, p. 282, Theorem 6.2.6] {\it
 Let $A_1,...,A_d  \in GL(s,\ZZ)$ be commuting partially hyperbolic
matrices with determinants $ w_1,..., w_d$, $p^{(k)}$ the product of the first $k$
primes numbers relatively prime to $ w_1,..., w_d$. If $ \bz,\tilde{\bz} \in \ZZ^s$  and there are
 $d$ sequences $ \{j_i^{(k)}, 1 \le i \le d  \}$ of integers such that
\begin{equation}\nonumber
   A_1^{j_1^{(n)}}  \cdots   A_d^{j_d^{(k)}} \tilde{\bz} \equiv \bz (\mod  p^{(k)}), \qquad k=1,2,...,
\end{equation}
then there exists a vector $(j_1^{(0)},...,j_d^{(0)}) \in \ZZ^s$ such that}
\begin{equation}\label{3.10}
   A_1^{j_1^{(0)}}  \cdots   A_d^{j_d^{(0)}} \tilde{\bz} = \bz .
\end{equation}

Let
\begin{equation}\label{3.12}
     \mu =   \begin{cases}
    1,  & \; {\rm if}  \;  s  \;{\rm is \;odd},\\
    2,   & \; {\rm if}  \;  s  \;{\rm is \;even, \;and\;} \nexists \bepsilon \; {\rm with} \; N_{\cK/Q} (\bepsilon) =-1, \\
    3,   & \; {\rm if}  \;  s  \;{\rm is \;even, \;and\;} \exists \bepsilon \; {\rm with} \; N_{\cK/Q} (\bepsilon) =-1. \\
  \end{cases}
\end{equation}

Let $\mu \in \{ 1,2 \}$.  By [BS, p. 117], we see that there exist units  $\bepsilon_{i} \in \cU_{\fO} $, with $N_{\cK/Q} (\bepsilon_i) =1$,  $i=1,...,s-1$, such that 
 every  $\bepsilon \in \cU_{\fO} $ can be uniquely represented as follows
\begin{equation}\label{3.15}
   \bepsilon = (-1)^a\bepsilon_{1}^{k_1} ...\bepsilon_{s-1}^{k_{s-1}},  \quad {\rm with} \quad (k_1,...,k_{s-1}) \in \ZZ^{s-1}, \; a \in \{0,1\}. 
\end{equation} 
Let $\mu=3$.  By [BS, p.~117],  there exist units  $\bepsilon_{i} \in \cU_{\fO} $, with $N_{\cK/Q} (\bepsilon_i) =1$,  $i=1,...,s-1$ and $N_{\cK/Q} (\bepsilon_0) =-1$, such that 
 every  $\bepsilon \in \cU_{\fO} $ can be uniquely represented as follows
\begin{equation}\label{3.16}
   \bepsilon = (-1)^{a_1} \bepsilon_{0}^{a_2}\bepsilon_{1}^{k_1} ...\bepsilon_{s-1}^{k_{s-1}},  \quad {\rm with} \quad (k_1,...,k_{s-1}) \in \ZZ^{s-1}, \; a_1,a_2 \in \{0,1\}.
\end{equation}

Consider the case $\mu=1$.  Let $ I_i ={\rm diag}((\sigma_j(\bepsilon_i))_{1 \leq j \leq s}), \; i=1,...,s-1   $, $\Gamma_{\fO}=\sigma(\fO)$,
 $\bff_1,...,\bff_s$ be a basis of $\Gamma_{\fO}$, $\be_i=(0,...,1,...,0) \in \ZZ^s$,  $i=1,...,s$ a basis of $\ZZ^s$. Let $Y$ be the $s \times s$ matrix with $\be_i Y = \bff_i$, 
  $i=1,...,s$.  We have $\ZZ^s Y = \Gamma_{\fO}$.  Let $A_i = Y I_i  Y^{-1}$,  $i=1,...,s-1   $. We see $\ZZ^s A_i =\ZZ^s  $ $(i=1,...,s-1)   $.
  Hence, $A_i$ is the integer matrix with $\det A_i = \det I_i =1 $  $(i=1,...,s-1 )  $.
  
 Let $\tilde{\bz} =(1,...,1)$ and $\bz =-\tilde{\bz} $. Let $h>1$, and let $A_s=p_2I$, where $I$ is the identity matrix.  
  Taking into account that $(\bepsilon_{1}^{k_1} ...\bepsilon_{s-1}^{k_{s-1}}p_2^{k_{s}})_j =1$ for some $j \in [1,s]$ if and only if $k_1=...=k_{s} =0$, we get that 
  $ A_{1},...,A_{s}$  are commuting  partially hyperbolic matrices.  By Definition 4, $-1$ is not the eigenvalue of $ A_{1}^{k_1} ...A_{s}^{k_{s}}$, and  $ \tilde{\bz}  A_{1}^{k_1} ...A_{s}^{k_{s}}  \neq  \bz  $ for  all $(k_1,...,k_{s}) \in \ZZ^{s}$.  Applying Theorem D with $d=s$,
  we have that there exists an integer $i \geq 1$  such that $(p_2,p^{(i)})=1$,
\begin{equation} \nonumber
   \tilde{\bz}  A_{1}^{k_1} ...A_{s-1}^{k_{s-1}}  \not\equiv  \bz  (\mod p^{(i)})  \quad {\rm for \; all} \quad (k_1,...,k_{s-1}) \in \ZZ^{s-1}, 
\end{equation}   
and
\begin{equation}\label{3.20}
     (\bepsilon_{1}^{k_1} ...\bepsilon_{s-1}^{k_{s-1}})_j \not\equiv -1 (\mod p^{(i)})  \quad {\rm for \; all} \quad (k_1,...,k_{s-1}) \in \ZZ^{s-1}, \quad j\in [1,s].
\end{equation} 
We denote this $p^{(i)}$ by $p_3$. We have $(p_2,p_3)=1$. If $h=1$, then we apply Theorem~D with $d=s-1$.

 Let $\fp_3 =p_3\fO$ and $\PP = \fO/\fp_3$. Denote the projection
 map $\fO \to \PP$  by $\pi_1$.   Let $\fO^*$ be the set of all integers of $\fO$ which are relatively prime to $\fp_3$, 
 $\PP^*=\pi_1(\fO^*)$,
\begin{equation}\nonumber
    \cE_j  = \{ v \in \PP^* \; | \; \exists \; (k_1,...,k_{s-1}) \in \ZZ^{s-1} \; {\rm with}  \; v \equiv  (-1)^j \bepsilon_{1}^{k_1} ...\bepsilon_{s-1}^{k_{s-1}}  (\mod  \fp_3) \},
\end{equation} 
where  $j =0,1$,  and   $\cE=\cE_0 \cup \cE_1$. By (\ref{3.20}),  $\cE_0 \cap \cE_1 =\emptyset$. Let
\begin{equation}\label{3.24}
    \chi_{1,p_3}(v) = (-1)^j,    \quad {\rm for} \quad  v \in \cE_j, \; j=0,1.
\end{equation} 
We see that $\chi_{1,p_3}$ is the character on  group $\cE$. We need the following known assertion (see e.g. [Is, p. 63], 
[Ko, p.~446, Ch.8, $\mathsection 2$, Ex.4]) :\\

 {\bf Lemma B.} {\it  Let $\dot{G}$  be a finite abelian group, $\dot{H}$ is a subgroup of $\dot{G}$, and $\chi_{\dot{H}}$ is a  character of $\dot{H}$. Then there exists a character $\chi_{\dot{G}}$
  of $\dot{G}$ such that $\chi_{\dot{H}}(h) = \chi_{\dot{G}}(h)$ for all $h \in \dot{H}$. }

Applying Lemma B, we can extend  the character $\chi_{1,p_3}$ to a character $\chi_{2,p_3}$ of group $\PP^*$. Now  we
extend $\chi_{2,p_3}$ to a character $\chi_{3,p_3}$  of group $\fO^*$ by setting
\begin{equation}\label{3.26}
    \chi_{3,p_3}(v) =   \chi_{2,p_3}(\pi_1(v)) \quad {\rm for} \quad  v \in \fO^*.
\end{equation} 
Let
\begin{equation}\nonumber
    \chi_{4,p_3}(v) =   \chi_{3,p_3}(v) \chi_{\infty}(v) \quad {\rm with} \quad  \chi_{\infty}(v) =\Nm(v)/|\Nm(v)| ,    
\end{equation} 
 for $ v \in \fO^*$, and let
\begin{equation}\label{3.29}
    \chi_{5,p_3}((v)) =   \chi_{4,p_3}(v).  
\end{equation}  
We need to prove that the right hand side of (\ref{3.29}) does not depend on units $\bepsilon \in \cU_{\fO}$. Let $\bepsilon=\bepsilon_{1}^{k_1} ...\bepsilon_{s-1}^{k_{s-1}}$. 
By (\ref{3.15}), (\ref{3.24}), and (\ref{3.26}), we have $\chi_{3,p_3}(\bepsilon) =1$, $\Nm(\bepsilon)=1$, and $\chi_{\infty}(\bepsilon)=1$. Therefore
\begin{equation}\nonumber
    \chi_{4,p_3}(v\bepsilon) =   \chi_{3,p_3}(v\bepsilon) \chi_{\infty}(v\bepsilon) = \chi_{3,p_3}(v) \chi_{3,p_3}(\bepsilon) \chi_{\infty}(v)\chi_{\infty}(\bepsilon) = \chi_{3,p_3}(v) \chi_{\infty}(v).
\end{equation} 
Now let $\bepsilon =-1$. Bearing in mind that $\chi_{3,p_3}(-1) =-1$, $\Nm(-1)=-1$, and $\chi_{\infty}(-1)=-1$, we obtain  $ \chi_{4,p_3}(-1)=1$. 
Hence, definition (\ref{3.29}) is correct. Let $\cI^{\fp_3}$ be the group
of all principal ideals of $\cK$ which are relatively prime to $\fp_3$. Let 
\begin{equation}\nonumber
    \chi_{6,p_3}((v_1/v_2)) =   \chi_{5,p_3}((v_1))/\chi_{5,p_3}((v_2)) \quad  {\rm for}  \quad v_1,v_2  \in \fO^*.
\end{equation} 
Let  $\cP^{\fp_3}$ is the group of
fractional principal ideals $(a)$ such that
 $a \equiv 1 \mod \fp_3$ and $\sigma_i(a) >0, \; i=1,...,s$.
Let  $\pi_2: \; \cI^{\fp_3} \to \cI^{\fp_3}/ \cP^{\fp_3}$ be the projection map. Bearing in mind that $\chi_{6,p_3}(\fa) = 1$ for $\fa \in  \cP^{\fp_3} $,
we define 
\begin{equation}\nonumber
    \chi_{7,p_3}(\pi_2(\fa)) = \chi_{6,p_3}(\fa)    \quad  {\rm for}  \quad \fa  \in  \cI^{\fp_3} .
\end{equation} 
By [Na, p. 94, Lemma 3.3],  $\cJ^{\fp_3}/ \cP^{\fp_3}$ is the finite abelian group.
Applying Lemma~B, we extend the character $\chi_{7,p_3}$ to  a character $\chi_{8,p_3}$  of  group $\cJ^{\fp_3}/ \cP^{\fp_3}$.
 We have  $\chi_{8,p_3}(\fa) = 1$ for $\fa \in  \cP^{\fp_3} $, and we set $ \chi_{9,p_3}(\fa) = \chi_{8,p_3}(\pi_3(\fa)) $, where $\pi_3$ is the proection map $\cJ^{\fp_3}    \to \cJ^{\fp_3}/ \cP^{\fp_3}$.   
It is easy to verify
\begin{equation}\nonumber
    \chi_{9,p_3}((v)) = \chi_{8,p_3}(\pi_3((v))) = \chi_{7,p_3}(\pi_3((v)))=\chi_{7,p_3}(\pi_2((v))) 
\end{equation} 
\begin{equation}\nonumber
 =\chi_{6,p_3}((v))=\chi_{4,p_3}(v) =\chi_{3,p_3}(v) \chi_{\infty}(v)       
\end{equation} 
  for  \quad $\fa  \in  \cI^{\fp_3} $.  Thus we have constructed a nontrivial Hecke character.
  
Case $\mu=2$. We repeat the construction of the case $\mu=1$, taking  $p_3=1$ and $ \chi_{4,p_3}((v)) =  \Nm(v)/|\Nm(v)|$. 

Case $\mu=3$.  Similarly to the case $\mu=1$, we have that there exists $i>0$ with
\begin{equation}\label{3.30}
     \bepsilon_{1}^{k_1} ...\bepsilon_{s-1}^{k_{s-1}} \not\equiv \bepsilon_0 (\mod p^{(i)}),  \quad {\rm for \; all} \quad (k_1,...,k_{s-1}) \in \ZZ^{s-1}.
\end{equation} 
We denote this $p^{(i)}$ by $p_3$. Let 
\begin{equation}\nonumber
    \cE_j  = \{ v \in \cP^* \; | \; \exists \; (k_1,...,k_{s-1}) \in \ZZ^{s-1} \; {\rm with}  \; v \equiv  \bepsilon_0^j \bepsilon_{1}^{k_1} ...\bepsilon_{s-1}^{k_{s-1}}  (\mod p_3 \fO) \},
\end{equation} 
where  $j =0,1$,  and  $\cE=\cE_0 \cup \cE_1$.  By (\ref{3.30}),  $\cE_0 \cap \cE_1 =\emptyset$.   Let
\begin{equation}\nonumber
    \chi_{2,p_3}(v) = (-1)^j    \quad {\rm for} \quad  v \in \cE_j, \; j=0,1.
\end{equation} 
Next, we repeat the construction of the case $\mu=1$, and we verify the correction of definition (\ref{3.29}).
Thus, we have proved the following lemma:\\

{\bf Lemma 4.} {\it Let $\mu \in \{1,2,3\}$. There exists $p_3 =p_3(\mu) \geq 1$, $(p_2,p_3)=1,$   a nontrivial Hecke character  $\dot{\chi}_{p_3}$, and a character $\ddot{\chi}_{p_3}$ on group $(\fO/p_3\fO)^*$ such that 
\begin{equation}\nonumber
    \dot{\chi}_{p_3} ((v)) =\tilde{\chi}_{p_3}(v)  ,  \quad {\rm with } \quad  \tilde{\chi}_{p_3}( v)  = \ddot{\chi}_{p_3}(v) \Nm(v)/|\Nm(v)|  ,  
\end{equation} 
for  $v \in \fO^*$, and $ \ddot{\chi}_{p_3}(v) =0 $ for $( v, p_3\fO)  \neq 1$}.
\\
\\

{\bf 2.4. Non-vanishing of $L$-functions.} 
With every Hecke character $\chi \; \mod \fm$, we associate its $L$-function
\begin{equation}\nonumber
   L(s, \chi)= \sum_{\fa }  \frac{\chi(\fa)}{\fN (\fa)^s} , 
\end{equation}
where $\fa$ varies over the integral ideals of $\cK$,   and we put $\chi(\fa) = 0$ whenever $(\fa, \fm) \neq  1$.   
\\

{\bf Theorem C.} [La, p. 313, Theorem 2]. Let $\chi$ be a nontrivial Hecke character. Then
\begin{equation}\nonumber
     L(1, \chi) \neq  0.  
\end{equation}
\\

{\bf Theorem D.} [RM, p. 128,  Theorem 10.1.4]  Let $(a_k)_{k \geq 1}$ be a sequence of complex numbers, and let
  $ \sum_{k<x} a_k = O(x^{\delta})$, for some $\delta > 0$. Then 
\begin{equation}\label{4.2}
       \sum_{n \geq 1} a_n/n^s
\end{equation}
converges for $ \Re(s) > \delta$. \\

{\bf Theorem E.} [Na, p. 464, Proposition I]
 If the series (\ref{4.2}) converges at a point $s_0$, then it converges
also in the open half-plane $ \Re s > \Re s_0 $, the convergence being uniform in
every angle $ \arg(s-s0) < c < \pi/2$. Thus (\ref{4.2}) defines a function regular in $ \Re s > \Re s_0 $. \\
 
Let $\bff_1,...,\bff_s$ be a basis of $\Gamma_{\fO}$, and let   $\bff_1^{\bot},...,\bff_s^{\bot}$  be a dual basis (i.e.   $ \langle \bff_i,\bff_i^{\bot}\rangle=1$, $ \langle \bff_i,\bff_j^{\bot} \rangle =0,$  $1 \leq i,j \leq s, \; i \neq j$   ).
Let 
\begin{equation} \label{4.2a}
   \Lambda_w =   \{ a_1 \bff_1^{\bot}+ \cdots + a_s \bff_s^{\bot} \;  |  \;   0 \leq a_i  \leq w-1, \;  i=1,...,s  \}, 
\end{equation}
 and $  \Lambda_w^* =  \{ \bb \in \Lambda_w\; | \; (w,\bb) =1 \} $. \\

{\bf Lemma 5.} {\it With notations as above,
\begin{equation}\label{4.4}
  \rho(M,j):=   \sum_{ \bgamma \in \Gamma_{j} \cap \cF, \;  |\Nm (\bgamma)|  < M/2} \tilde{\chi}_{ p_3}(\bgamma) =  O(M^{1-1/s}), \quad  \quad j \in [1,h],
\end{equation}  
and
\begin{equation}\label{4.6}
  \sum_{  \fN(\fa)  < M/2} \dot{\chi}_{ p_3} (\fa) =  O(M^{1-1/s})
\end{equation} 
for $M \to \infty$, where $\fa$ varies over the integral ideals of $\cK$. } \\

 {\bf Proof.} By Lemma 4, we have
\begin{equation}\nonumber
       \rho(M,j) = \sum_{\ba  \in \Lambda^{*}_{p_3} } \ddot{\chi}_{ p_3}(\ba) \sum_{\varsigma_i \in \{-1, +1\},\; i=1,...,s}      \varsigma_1 \cdots \varsigma_s  \dot{\rho}(\ba, \bvarsigma,j),
\end{equation} 
where
\begin{equation}\nonumber
  \dot{\rho}(\ba, \bvarsigma,j)  =   \sum_{ \substack{\bgamma \in \Gamma_{j} \cap \cF, \; \bgamma \equiv \ba \mod  p_3,  \\  |\Nm (\bgamma)|  < M/2, \;\sgn(\gamma_i) =\varsigma_i, \;i=1,...,s }} 1  .
\end{equation} 
Using Lemma 3  with $M_1=M/2$ and   $w=p_3$ , we get
\begin{equation}\nonumber
      \dot{\rho}(\ba, \bvarsigma,j)   =  \sum_{ \substack{\bgamma \in (p_3\Gamma_{j} +\ba) \cap \cF, \;|\Nm (\bgamma)|  < M/2 \\   \sgn(\gamma_i) =\varsigma_i, \;i=1,...,s }} 1 
   =c_{6,j} M/p_3^s + O(M^{1-1/s}).
\end{equation}
Therefore
\begin{equation}\nonumber
       \rho(M,j) = \sum_{\ba \in \Lambda^{*}_{p_3} } \ddot{\chi}_{ p_3}(\ba) \sum_{\varsigma_i \in \{-1, +1\},\; i=1,...,s}      \varsigma_1 \cdots \varsigma_s 
   (c_{6,j} M /p_3^s + O(M^{1-1/s}))   = O(M^{1-1/s}).
\end{equation}
Hence,  the assertion (\ref{4.4}) is proved.  The assertion (\ref{4.6})  can be proved similarly
(see also [CaFr, p. 210, Theorem 1], [Mu, p. 142, and p.144, Theorem 11.1.5]). 
  \quad  \qed \\

{\bf Lemma 6.} {\it There exists $M_0>0$, $i_0 \in [1,h]$, and $c_{7}>0$, such that}
\begin{equation}\nonumber
   |\rho_0(M,i_0)|   \geq c_{7}  \quad {\rm  for} \quad M>M_0, \quad {\rm  with} \quad   \rho_0(M,i)=\sum_{ \bgamma \in \Gamma_{i} \cap \cF, \;  |\Nm (\bgamma)|  < M/2}  \frac{\tilde{\chi}_{p_3}(\bgamma)}{|\Nm (\bgamma)|}.
\end{equation} 
 {\bf Proof.} Let $ cl(\cK) = \{C_1,..., C_h \}$, $\fa_i \in C_i$  be an integral ideal, $i=1,...,s$,  
 and let $C_1$ be the class of principal ideals. Consider the inverse ideal class $C_i^{-1}$. We set $\dot{\fa}_i =\{\fa_1,...,\fa_h \} \cap C_i^{-1}$.
Then  for any $\fa \in C_i$ the product $\fa \dot{\fa}_i$ will be a principal ideal: $\fa \dot{\fa}_i = (\alpha)$, $(\alpha \in \cK)$. 
By [BS, p. 310], we have that the mapping $\fa \to (\alpha)$ establishes a one to one correspondence between integral ideal $\fa$ of the class
 $C_i$ and principal ideals divisible by $\dot{\fa}_i$. Let
\begin{equation}\nonumber
  \rho_1(M) = \sum_{\fN (\fa)  < M/2 } \dot{\chi}_{p_3}(\fa)/\fN(\fa) .
\end{equation}
 Similarly to [BS, p. 311], we get 
\begin{equation}\nonumber
  \rho_1(M) = \sum_{1 \leq i \leq h}  \;\;  \sum_{\fa \in C_i, \fN (\fa)  < M/2 }  \frac{\dot{\chi}_{p_3}(\fa)}{\fN(\fa)} 
  = \sum_{1 \leq i \leq h}  \;\;  \sum_{\substack{ \fa \in C_1, \fN (\fa/ \dot{\fa}_i))  < M/2 \\ \fa \equiv 0 \mod \dot{\fa}_i }} 
    \frac{\dot{\chi}_{p_3}(\fa/ \dot{\fa}_i)}{\fN(\fa/ \dot{\fa}_i))} .
\end{equation}
Let
\begin{equation}\nonumber
  \rho_2(M,i) =   \sum_{\substack{ \fa \in C_1, \; \fN (\fa)  < M/2 \\ \fa \equiv 0 \mod \dot{\fa}_i }} \dot{\chi}_{p_3}(\fa)/\fN(\fa) .
\end{equation}
We see
\begin{equation}\label{4.7}
  \rho_1(M) = \sum_{1 \leq i \leq h}  \frac{ \dot{\chi}_{p_3}(1/  \dot{\fa}_i)}{\fN(1/ \dot{\fa}_i)}  \rho_2(M\fN( \dot{\fa}_i) ,i).
\end{equation}
By Lemma 4, we obtain $\tilde{\chi}_{p_3}(\bgamma) /|\Nm (\bgamma)| = \dot{\chi}_{p_3}((\bgamma)) /\fN ((\bgamma))$. 
Using  Theorem B, we get $\rho_0(M,i) = \rho_2(M,i)$.
From (\ref{4.6}), Theorem C,   Theorem D, and  Theorem~E, we derive    $ \rho_1(M)   \stackrel{M \rightarrow\infty}{\longrightarrow} L(1, \dot{\chi}_{p_3}) \neq 0 $. By (\ref{4.4}) and Theorem D, we obtain that  there exists a complex number $\rho_i$, such that 
$ \rho_0(M,i)   \stackrel{M \to \infty}{\longrightarrow} \rho_i $, $i=1,...,h$.
Hence,  there exists $M_0>0$, such that 
\begin{equation}\label{4.8}
 |L(1, \dot{\chi}_{p_3})|/2  \leq  |\rho_1(M)| , \quad {\rm and} \quad |\rho_i -\rho_2(M,i)| \leq |L(1, \dot{\chi}_{p_3})| (8 \beta)^{-1},
\end{equation}
with  $\beta=  \sum_{1 \leq i \leq h}   \fN(\dot{\fa}_i)$,   for $M  \geq M_0$.
 Let $\rho = \max_{1 \leq i \leq h} |\rho_i| = |\rho_{i_0}|  $.\\ 
Using  (\ref{4.7}), we have 
\begin{equation}\nonumber
    |L(1, \dot{\chi}_{p_3})|/2  \leq |\rho_1(M)| 
   \leq \rho \beta  + \Big| \sum_{1 \leq i \leq h} \frac{ \dot{\chi}_{p_3}(1/  \dot{\fa}_i)}{\fN(1/ \dot{\fa}_i)} 
   (\rho_i -  \rho_2(M \fN( \dot{\fa}_i)  ,i) ) \Big|     
\end{equation}
\begin{equation}\nonumber
   \leq \rho \beta +  |L(1, \dot{\chi}_{p_3})|/8 \quad {\rm for}\quad M>M_0.
\end{equation}
By (\ref{4.8}), we get for $M>M_0$
\begin{equation}\nonumber
       \rho  \geq |L(1, \dot{\chi}_{p_3})|(4 \beta)^{-1},   \quad  {\rm and} \quad   |\rho_0(M,i_0)|   =    |\rho_2(M,i_0)| \geq |L(1, \dot{\chi}_{p_3})|(8 \beta)^{-1}.
\end{equation}
Therefore, Lemma 6 is proved.~\qed \\

{\bf Lemma 7.} {\it There exists $M_2>0$, such that}
\begin{equation}\nonumber
   |\vartheta |  \geq c_{7}/2  \quad {\rm  for} \quad M>M_2, \quad {\rm  where} \quad
   \vartheta = \sum_{ \bgamma \in \Gamma_{i_0} \cap \cF}  \frac{\ddot{\chi}_{p_3}(\bgamma) \eta_M(\bgamma)}{\Nm (\bgamma)} .
\end{equation} 
 {\bf Proof.} Let $\dot{\eta}_M(k)=1-\eta(2|k|/M)$, 
\begin{equation}\nonumber
    \vartheta_1= \sum_{\substack{ \bgamma \in \Gamma_{i_0} \cap \cF \\  |\Nm (\bgamma)| < M/2}}  \frac{\tilde{\chi}_{p_3}(\bgamma) }{|\Nm (\bgamma)|} ,
      \quad {\rm  and} \quad
    \vartheta_2 = \sum_{ \substack{\bgamma \in \Gamma_{i_0} \cap \cF \\  M/2 \leq |\Nm (\bgamma)|  \leq M}}  \frac{\tilde{\chi}_{p_3}(\bgamma) \dot{\eta}_M(\Nm(\bgamma))}{|\Nm (\bgamma)|}.
\end{equation}   
From (\ref{2.35a}), we get $\eta_M(\bgamma) =\dot{\eta}_M(\Nm(\bgamma))$, $\eta_M(\bgamma) =1$  for $|\Nm(\bgamma)| \leq M/2$, and  $\eta_M(\bgamma) =0$  for $|\Nm(\bgamma)| \geq M$. Using Lemma 4, we derive
\begin{equation}\label{4.30}
  \vartheta = \sum_{ \bgamma \in \Gamma_{i_0} \cap \cF, \;  |\Nm (\bgamma)|  \leq M}  \frac{\tilde{\chi}_{p_3}(\bgamma) \dot{\eta}_M(\Nm(\bgamma))}{|\Nm (\bgamma)|} , \quad {\rm  and} \quad  \vartheta =\vartheta_1 + \vartheta_2.
\end{equation}  
Bearing in mind that $\Nm(\bgamma) \in \ZZ$ and $\Nm(\bgamma) \neq 0$, we have
\begin{equation} \nonumber
       \vartheta_2 = \sum_{  M/2 \leq \dot{n} \leq M}  \frac{a_{\dot{n}} \dot{\eta}_M(k)}{k}  , \quad {\rm  with} \quad  a_{\dot{n}} = \sum_{ \bgamma \in \Gamma_{i_0} \cap \cF, \;  |\Nm (\bgamma)|  =\dot{n}}  \tilde{\chi}_{p_3}(\bgamma) .
\end{equation}  
Applying Abel' transformation
\begin{equation}\nonumber
       \sum_{m < k \leq \dot{n}}  g_k f_k  =g_{\dot{n}} F_{\dot{n}}- \sum_{m < k \leq \dot{n}-1}  (g_{k+1} -g_k)F_k , \quad {\rm  where} \quad F_k= \sum_{m < i \leq k}  f_i,
\end{equation}  
with $f_k = a_k, \; g_k = \dot{\eta}_M(k)/k$ and       $F_k=  \sum_{ \bgamma \in \Gamma_{i_0} \cap \cF, \;  M/2 -0.1 < |\Nm (\bgamma)|  \leq k }  \tilde{\chi}_{p_3}(\bgamma)$, \\we obtain
\begin{equation}\label{4.32}
       \vartheta_2 = \dot{\eta}_M(M) F_{M}/M - \sum_{M/2 -0.1< k \leq M-1}  (\dot{\eta}_M(k+1)/(k+1) -  \dot{\eta}_M(k)/k) F_k .
\end{equation}  
Bearing in mind that $ 0 \leq \dot{\eta}_M(x) \leq 1  $ and $\eta^{'}(x) =O(1)$, for $|x| \leq 2$, we get 
\begin{equation} \nonumber
    |\dot{\eta}_M(k+1)/(k+1) - \dot{\eta}_M(k)/k)| \leq    |\dot{\eta}_M(k+1)/(k+1) -  \dot{\eta}_M(k+1)/k)|  
\end{equation}  
\begin{equation} \nonumber
      +  |(\dot{\eta}_M(k+1)-\dot{\eta}_M(k))/k|  \leq 1/k^2 + 2(kM)^{-1} \sup_{x \in [0,2]} |\eta^{'}(x)| = O(k^{-2}).
\end{equation} 
Taking into account that $F_k = O(M^{1-1/s}) $ (see (\ref{4.4})), we have from 
 (\ref{4.32}) that $\vartheta_2 =O(M^{-1/s})$. Using Lemma 6 and (\ref{4.30}), we obtain the assertion of  Lemma~7. \\ \hfill \qed  \\
\\

{\bf 2.5. The  lower bound estimate for $\bE(\cA(\bx,M))$.} 
 Let $n=s^{-1} \log_2 N$ with $N=N_1 \cdots N_s$,  $\tau =N^{-2}$,  $M=[\sqrt{n}]$,  and
\begin{equation}\nonumber
    G_0=\{ \bgamma \in \Gamma^\bot \; |  \;  |\Nm(\bgamma)|  > M  \}, 
\end{equation}
\begin{equation}\nonumber
      G_1=\{ \bgamma \in \Gamma^\bot \; |  \;  |\Nm(\bgamma)|  \leq M , \; \max_i |\gamma_i|  \geq 1/\tau^2  \},
\end{equation}
\begin{equation}\label{5.0}
    G_2=\{ \bgamma \in \Gamma^\bot \; |  \; |\Nm(\bgamma)|  \leq M , \; 1/\tau^2 > \max_i |\gamma_i| \geq n/\tau  \},
\end{equation}
\begin{equation}\nonumber
   G_3=\{ \bgamma \in \Gamma^\bot \; |  \; |\Nm(\bgamma)|  \leq M , \;  n/\tau   > \max_i |\gamma_i| \geq n^{-s}/ \tau \},
\end{equation}
\begin{equation}\nonumber
     G_4=\{ \bgamma \in \Gamma^\bot  \setminus {\bs}\; |  \; |\Nm(\bgamma)|  \leq M , \; \max_i|\gamma_i| <n^{-s}\tau^{-1},\;\;       n^{-s}  >  N^{1/s}\min_i|\gamma_i| \},
\end{equation}
\begin{equation}\nonumber
   G_5=\{ \bgamma \in \Gamma^\bot  |  \; |\Nm(\bgamma)|  \leq M , \;   \max_i|\gamma_i| <n^{-s}\tau^{-1} , \;   N^{1/s} \min_i|\gamma_i| \in [  n^{-s},n^{s}]   \},
\end{equation}
\begin{equation}\nonumber
   G_6=\{ \bgamma \in \Gamma^\bot \; |  \; |\Nm(\bgamma)|  \leq M , \;\;  \max_i|\gamma_i| <n^{-s}\tau^{-1}, \;\quad \; N^{1/s}\min_i|\gamma_i| >  n^{s} \}.
\end{equation}
We see that
\begin{equation}\nonumber
              \Gamma^\bot  \setminus {\bs} = G_0 \cup \cdots \cup G_6, \quad {\rm and} \quad G_i \cap G_j =\emptyset, \; {\rm for}\; i\neq j.
\end{equation} 

 Let $p=p_1 p_2 p_3$, $\bb \in \Delta_p$. By (\ref{2.35a}) and (\ref{2.35}), we have 
\begin{equation}\label{5.2}
    \cA(\bb/p,M)   = \sum_{0 \leq i \leq 6} \cA_i(\bb/p,M), \quad {\rm and} \quad  \cA_0(\bb/p,M)=0,
\end{equation} 
where
\begin{equation}\label{5.2a}
    \cA_i(\bb/p,M)   = \sum_{\bgamma\in G_i } 
          \prod_{i=1}^s \sin(\pi \theta_i N_i \gamma_i) \frac{    \eta_M(\bgamma) \widehat{\Omega} (\tau \bgamma)e(\langle\bgamma,\bb/p\rangle + \dot{x})}{\Nm(\bgamma)},
\end{equation}
with $\dot{x} = \sum_{1 \leq i  \leq s}  \theta_i N_i  \gamma_i/2  $.

We will use the following simple decomposition (see notations from   $\mathsection 2.2$ and  (\ref{3.12}) - (\ref{3.16})):
\begin{equation} \nonumber
   G_i =   \bigcup_{1 \leq j \leq M} \;\;\bigcup_{ \bgamma_0 \in \Gamma^\bot \cap \cF, |\Nm(\bgamma_0)| \in (j-1,j ]  }  \;\; \bigcup_{a_1,a_2 =0,1}
      \Big\{ \bgamma \in G_i \;\;|  
\end{equation}
\begin{equation} \label{5.3}
    \bgamma=  \bgamma_0  (-1)^{a_1}\bepsilon_{0}^{a_2}\bepsilon^{\bk} , \;   \bk \in \ZZ^{s-1} \Big\}, \quad     i \in [1,6] ,           
\end{equation} 
 where $\bk =(k_1,...,k_{s-1})$, $\bepsilon^{\bk} = \bepsilon_1^{k_1} \cdots \bepsilon_{s-1}^{k_{s-1}} $, and $\bepsilon_0=1$ for $\mu=1,2$.\\

{\bf Lemma 8.} {\it With notations as above} 
\begin{equation} \nonumber
    \cA_i(\bb/p, M) =O( n^{s-3/2}\ln n), \quad {\rm where } \quad  M= [\sqrt{n}]  \quad {\rm and } \quad   i \in [1,5].
\end{equation} 
 {\bf Proof.} By (\ref{5.2a}), we have
\begin{equation}\label{5.6}
    |\cA_i(\bb/p,M)| \leq   \sum_{\bgamma\in G_i}
            \prod_{1 \leq j \leq s} \big|\sin(\pi \theta_j N_j\gamma_j)   \widehat{\Omega} (\tau  \bgamma)/ \Nm(\bgamma) \big|.
\end{equation}

 {\it Case} $i=1.$  Applying (\ref{3.2}), we obtain $  \#\{\bgamma \in \Gamma^{\bot} \; : \; j \leq  | \bgamma | \leq j+1\}  =O(j^{s-1}) $.
By (\ref{2.12}) we get $  \widehat{\Omega} (\tau \bgamma) =O( (\tau |\bgamma|)^{-2s})  $ for $\bgamma \in G_1$.
From (\ref{5.6}) and (\ref{5.0}), we have 
\begin{equation}\nonumber
    \cA_1(\bb/p,M) =O\Big( \sum_{\bgamma\in \Gamma^\bot, \max_{i \in [1,s]}|\gamma_i| \geq 1/\tau^2 }
            \tau^{-2s}(\max_{i \in [1,s]}|\gamma_i|)^{-2s}  \Big)
\end{equation}
\begin{equation} \nonumber
   = O\Big( \sum_{j \geq\tau^{-2}}  \sum_{  \substack{ \bgamma\in \Gamma^\bot \\ \max_{i}|\gamma_i| \in [j,j+1) }}
            \tau^{-2s}(\max_{i \in [1,s]}|\gamma_i|)^{-2s}\Big) = O\Big( \sum_{j \geq \tau^{-2}} 
            \frac{\tau^{-2s}}{j^{s+1}}  \Big) =O(1).
\end{equation}

{\it Case} $i=2.$ By (\ref{2.12}) we obtain $  \widehat{\Omega} (\tau \bgamma) =O( n^{-2s})  $ for $\bgamma \in G_2$.
 By [BS, pp. 312, 322], the points of $\Gamma_\fO \cap \cF$ can be arranged in a sequence $\dot{\bgamma}^{(k)}$ so that 
\begin{equation}\label{5.7}
  |{\rm Nm}(\dot{\bgamma}^{(1)}) | \leq  |{\rm Nm}(\dot{\bgamma}^{(2)})| \leq ... \; {\rm and} \;
  c^{(1)} k   \leq |{\rm Nm}(\dot{\bgamma}^{(k)})|  \leq c^{(2)}  k, 
\end{equation} 
 $k=1,2,...$ for some $c^{(2)} > c^{(1)}>0$.  Let $ \bepsilon^{\bk}_{max} =\max_{1 \leq i \leq s}|(\bepsilon^{\bk})_i|$ and 
  $\bepsilon^{\bk}_{min} = \min_{1 \leq i \leq s}|(\bepsilon^{\bk})_i|$.
Using Lemma 2, we get 
\begin{equation} \label{5.9}
    \# \{ \bk \in \ZZ^{s-1} \; | \;  \bepsilon^{\bk}_{max}  \leq \tau^{-4}  \} = O(n^{s-1}), \quad {\rm where} \quad
    \tau=N^{-2}=e^{-2sn}.
\end{equation}
Applying (\ref{5.3}) - (\ref{5.9}), we have
\begin{equation} \nonumber
    \cA_2(\bb/p,M)  = O\Big(\sum_{1 \leq j \leq M} \sum_{ \bk \in \ZZ^{s-1}, \;\bepsilon^{\bk}_{max} \leq \tau^{-2} } 
           n^{-2s} \Big)  = O( M n^{-2s+s-1}) = O( 1) .
\end{equation}

{\it Case} $i=3.$ Using Lemma  2, we obtain 
\begin{equation}\label{5.9a}
   \#  \{ \bk \in \ZZ^{s-1} \; | \; \bepsilon^{\bk}_{max}
       \in [n^{-s-1}/\tau, n^{s+1}/\tau ]  \} 
\end{equation}
\begin{equation}\nonumber
    = c_4(\ln^{s-1} (n^{s+1}/\tau) - \ln^{s-1} (n^{-s-1}/\tau)  ) + O(n^{s-2}) 
\end{equation}
\begin{equation}\nonumber
        =  O\Big( |\ln^{s-1} \tau| \Big(\Big( 1 +\frac{(s+1)\ln n}{ |\ln \tau | }\Big)^{s-1} -   \Big( 1 -\frac{(s+1)\ln n}{ | \ln \tau | }
        \Big)^{s-1}\Big)\Big) = O(n^{s-2}\ln n).
\end{equation}
Applying (\ref{5.3}) - (\ref{5.9}), we get
\begin{equation} \nonumber
      \cA_3(\bb/p,M)  = O\Big(\sum_{1 \leq j \leq M} \;\;  \sum_{ \bk \in \ZZ^{s-1}, \bepsilon^{\bk}_{max}
       \in [n^{-s-1}/\tau, n^{s+1}/\tau ]}   \; 1 \Big)  = O( M n^{s-2} \ln n) .    
\end{equation}

{\it Case} $i=4.$ We see $\min_{1 \leq i \leq s} |\sin(\pi N_i \gamma_i)| =O(n^{-s})$ for $\bgamma \in G_4$.  Applying (\ref{5.3}) - (\ref{5.9}), we have
\begin{equation} \nonumber
    |\cA_4(\bb/p,M)| =    O\Big(  \sum_{1 \leq j \leq M} \sum_{ \bk \in \ZZ^{s-1}, \;\bepsilon^{\bk}_{max} \leq \tau^{-4}} 
          n^{-s} \Big)  = O( M n^{-2}).
\end{equation}

{\it Case} $i=5.$ 
Similarly to (\ref{5.9a}), we obtain from Lemma 2 that
\begin{equation}\nonumber
     \{ \bk \in \ZZ^{s-1} \; | \; \bepsilon^{\bk}_{min}   \in 
      [n^{-s-1}N^{-1/s} , n^{s+1}N^{-1/s} ]  \} 
= O(n^{s-2}\ln n).
\end{equation}
Therefore
\begin{equation} \nonumber
      \cA_3(\bb/p,M)  = O\Big(\sum_{1 \leq j \leq M}  \;\;  \sum_{ \bk \in \ZZ^{s-1}, \bepsilon^{\bk}_{min}
       \in [n^{-s-1}N^{-1/s} , n^{s+1}N^{-1/s} ] }   \; 1 \Big)  = O( M n^{s-2} \ln n) .    
\end{equation}
Hence, Lemma 8 is proved. \qed

Let $ \bvarsigma =(\varsigma_1,...,\varsigma_s)   $, $\b1=(1,1,...,1)$, and
\begin{equation}\label{5.18}
  \breve{ \cA}_6(\bb/p,M,\bvarsigma) =  \varsigma_1 \cdots \varsigma_s 
  (2\sqrt{-1})^{-s}  \sum_{\bgamma\in G_6 }
            \frac{\widehat{\Omega} (\tau  \bgamma) \eta_M(\bgamma)e(\langle \bgamma, \bb/p+\dot{\btheta}(\bvarsigma) \rangle)}{\Nm (\bgamma)}
\end{equation}
with $\dot{\btheta}(\bvarsigma) = (\dot{\theta}_1(\bvarsigma),...,\dot{\theta}_s(\bvarsigma))$ and $\dot{\theta}_i(\bvarsigma) = (1+\varsigma_i) \theta_iN_i/4 , \; i=1,...,s$.

By (\ref{5.2a}), we see
\begin{equation}\label{5.19}
    \cA_6(\bb/p,M) =  
   \sum_{ \bvarsigma \in \{1,-1\}^s}   \breve{ \cA}_6(\bb/p,M,\bvarsigma) .  
\end{equation}
\\

{\bf Lemma 9.} {\it With notations as above
\begin{equation}\nonumber
       \bE(\cA_6(\bb/p,M)) = \dot{\cA_6}(\bb/p,M,-\b1) + O(1),
\end{equation}
where}
\begin{equation}\label{5.20}
    \dot{\cA_i}(\bb/p,M,-\b1) = (-2\sqrt{-1})^{-s}  \sum_{\bgamma\in G_i }
            \frac{ \eta_M(\bgamma)e(\langle \bgamma,\bb/p\rangle)}{\Nm(\bgamma)},\quad i=1,2,...
\end{equation}
{\bf Proof.}  By  (\ref{5.18}) and (\ref{5.19}), we have
\begin{equation}\nonumber
   |{\bf E} (\cA_6(\bb/p,M)) -  \breve{ \cA}_6(\bb/p,M,-\b1) | =O\Big(  \sum_{ \substack{ \bvarsigma \in \{1,-1\}^s \\  \bvarsigma  \neq -\b1
    }}    
  \sum_{\bgamma\in G_6 }    \sum_{1 \leq i \leq s}
            \frac{      |{\bf E} (e(\varsigma_i \theta_iN_i \gamma_i/4))|}{|\Nm(\bgamma)|} \Big).
\end{equation}
 Bearing in mind that 
\begin{equation}\label{5.20a}
   {\bf E}(e(\theta_i  z )) =\frac{e(z)-1}{2\pi \sqrt{-1}  z}
\end{equation}      
and that $|N_i \gamma_i| \geq n^s/c_3$  for $\bgamma \in G_6$ (see (\ref{2.8}), and (\ref{5.0})),           we get
\begin{equation}\nonumber
      |{\bf E} (\cA_6(\bb/p,M)) -  \breve{ \cA}_6(\bb/p,M,-\b1) | =O\Big(
  \sum_{\bgamma\in G_6 }   
           n^{-s}|\Nm(\bgamma)|^{-1}\Big).
\end{equation} 
By (\ref{5.18}) and (\ref{5.20}), we obtain
\begin{equation}\nonumber
   |\breve{ \cA}_6(\bb/p,M,-\b1) - \dot{ \cA}_6(\bb/p,M,-\b1)| =O\Big(  
  \sum_{\bgamma\in G_6 }    
            \frac{      |\widehat{\Omega} (\tau  \bgamma) -1|}{|\Nm(\bgamma)|} \Big) .
\end{equation}
By (\ref{2.13})  and (\ref{5.0}), we see  $\widehat{\Omega} (\tau  \bgamma) =1 +O(n^{-s})$ for $\bgamma \in G_6$.
From  (\ref{5.0}),  (\ref{5.3}) and (\ref{5.9}), we have $\#G_6  =O(M n^{s-1})$. Hence
\begin{equation}\nonumber
       \bE(\cA_6(\bb/p,M)) - \dot{\cA}_6(\bb/p,M,-\b1) =  O\Big(\sum_{\bgamma \in G_6 }   
           n^{-s}|\Nm(\bgamma)|^{-1} \Big)=     O(1).
\end{equation}
Therefore, Lemma 9 is proved.  \qed

Let   
\begin{equation} \label{5.35}
 G_7 = \bigcup_{ \bgamma_0 \in \Gamma^{\bot} \cap \cF, |\Nm(\bgamma_0)| \leq M  }  \;\; \bigcup_{a_1,a_2 =0,1} \bigcup_{  \bk \in \cY_N}  T_{ \bgamma_0,a_1,a_2 ,\bk}  , 
\end{equation}
with
\begin{equation}\label{5.36}
    \cY_N  = \{ \bk \in \ZZ^{s-1} \; | \; \bepsilon^{\bk}_{min} \geq  N^{-1/s} \}, 
\end{equation}
and
\begin{equation}\nonumber
   T_{ \bgamma_0,a_1,a_2 ,\bk} =    \{ \bgamma \in \Gamma^\bot  \; | \; \bgamma=  \bgamma_0  (-1)^{a_1}\bepsilon_{0}^{a_2}\bepsilon^{\bk} \}.
\end{equation}
We note  that $ \# T_{ \bgamma_0,a_1,a_2 ,\bk} \leq 1  $ (may be $\bgamma_0  (-1)^{a_1}\bepsilon_{0}^{a_2}\bepsilon^{\bk} \notin  \Gamma^{\bot}$). \\

{\bf Lemma 10.} {\it With notations as above}
\begin{equation}\label{5.38}
       \bE(\cA(\bb/p,M)) = \dot{\cA_7}(\bb/p,M,-\b1) + O(n^{s-3/2} \ln n),   \quad {\rm where} \quad M=[\sqrt{n}] .
\end{equation}
{\bf Proof.}
By (\ref{5.20}), we have
\begin{equation}  \nonumber
      |\dot{\cA_6}(\bb/p,M,-\b1) -  \dot{\cA_7}(\bb/p,M,-\b1) | =O(\#(G_7 \setminus G_6) +  \#(G_6 \setminus G_7)) .
\end{equation} 
Consider $\bgamma \in G_6$ (see (\ref{5.0})). Bearing in mind that $ \min_{1 \leq i \leq s} |\gamma_i| \geq  n^{s} N^{-1/s}$, we get
\begin{equation}\nonumber
  |\gamma_i|  = |\Nm(\bgamma)| \prod_{[1,s] \ni j \neq i} |\gamma_j|^{-1} \leq  n^{-s(s-1)} N^{1+(s-1)/s}  < n^{-s}/\tau, \quad {\rm  with} \quad \tau =N^{-2}.
\end{equation}
Thus
\begin{equation}\nonumber
   G_6=\{ \bgamma \in \Gamma^\bot \; |  \; |\Nm(\bgamma)|  \leq M , \;\;  \;  N^{1/s}\min_i|\gamma_i| >  n^{s} \}.
\end{equation}
From   (\ref{5.35}), we obtain $G_7  \supseteq  G_6$.
Bearing in mind that $ |\Nm(\bgamma)| \geq 1$  for $\bgamma \in \Gamma^\bot  \setminus {\bs} $, we have that $G_6  \supseteq  G_5$, where
\begin{equation}\nonumber
 G_5 =  \bigcup_{ \bgamma_0 \in \Gamma^{\bot} \cap \cF, \;|\Nm(\bgamma_0)| \leq M  } \; \bigcup_{a_1,a_2 =0,1} \; \bigcup_{  \bk \in \dot{\cY}_N}  T_{ \bgamma_0,a_1,a_2 ,\bk}  , 
\end{equation}
with
\begin{equation}\label{5.36a}
    \dot{\cY}_N  = \{ \bk \in \ZZ^{s-1} \; | \;  N^{1/s} \bepsilon^{\bk}_{min} \geq n^{2s}\}.
\end{equation}
By  Lemma 3, we get  $ \# \{ \bgamma_0 \in \Gamma^{\bot} \cap \cF, |\Nm(\bgamma_0)|  \leq M \} =O(M).$
Therefore
\begin{equation} \nonumber
      |\dot{\cA}_6(\bb/p,M,-\b1) -  \dot{\cA}_7(\bb/p,M,-\b1) | =O(M \#( \cY_N \setminus \dot{\cY}_N) ) .
\end{equation}
Using Lemma  2, we obtain 
\begin{equation}\nonumber
    \#( \cY_N \setminus \dot{\cY}_N) = \{ \bk \in \ZZ^{s-1} \; | \; \bepsilon^{\bk}_{min}
       \in [N^{-1/s}, n^{2s} N^{-1/s}]  \} 
\end{equation}
\begin{equation}\nonumber
    = c_{5}\big(\ln^{s-1} ( N^{1/s}) - \ln^{s-1} (n^{-2s} N^{1/s})  \big) + O(n^{s-2}) 
\end{equation}
\begin{equation}\nonumber
        =  O\Big( \ln^{s-1} N \Big(\big(1 -   \big( 1 -\frac{2s^2\log_2 n}{  \ln N }   \big)^{s-1}\big)\Big)\Big) = O(n^{s-2}\ln n), \quad n=s^{-1}\log_2 N.
\end{equation}
Hence
\begin{equation} \nonumber
      |\dot{\cA_6}(\bb/p,M,-\b1) -  \dot{\cA_7}(\bb/p,M,-\b1) | =O(M n^{s-2} \ln n).
\end{equation}
Applying Lemma 8 and Lemma 9,    we get the assertion of Lemma 10. \qed  \\

Let 
\begin{equation}\nonumber
     \delta_w(\bgamma) =   \begin{cases}
    1,  & \; {\rm if}  \;   \bgamma  \in w\fO ,\\
     0,   & \; {\rm otherwise}  .
  \end{cases}
\end{equation} \\

{\bf Lemma 11.} {\it  Let $\bgamma \in \fO$, then
\begin{equation}\nonumber
    \frac{1}{w^s} \sum_{\by \in \Lambda_w}   e(<\bgamma, \by>/w)   = \delta_w(\gamma).
\end{equation} }
 {\bf Proof.} It easy to verify that
\begin{equation}\label{5.37}
    \frac{1}{v} \sum_{0 \leq k <w}   e(kb/w)   = \dot{\delta}_w(b),   \quad { \rm where} \quad   \dot{\delta}_w(b) =   \begin{cases}
    1,  & \; {\rm if}  \;  b \equiv 0 \mod w,\\
    0,   & \; {\rm otherwise}  .
  \end{cases}
\end{equation} 
Let $\bgamma=d_1 \bff_1+ \cdots + d_s \bff_s$, and $\by=a_1 \bff_1^{\bot}+ \cdots + a_s \bff_s^{\bot}$ (see (\ref{4.2a})). 
We have $\langle \bgamma, \by \rangle  = a_1d_1 + \cdots + a_s d_s$.
Bearing in mind that $ \bgamma \in w\fO  $ if and only if $d_i   \equiv 0 \mod w$ $(i=1,...,s)$, we obtain from (\ref{5.37}) the assertion of Lemma 11. \qed  \\

{\bf Lemma 12.} {\ There exist  $\bb \in \Lambda_{p}$, $c_8>0$ and $N_0>0$ such that
\begin{equation}\nonumber
    | \bE(\cA(\bb/p,M))|  > c_{8}   n^{s-1}     \quad  {\rm for } \quad N > N_0.
\end{equation} }
{\bf Proof.} We consider the case $\mu=1$. The proof for the cases $\mu=2,3$ is similar.
By (\ref{5.20}) and Lemma 11, we have 
\begin{equation}\nonumber
 \varrho:=  \frac{2^{2s}}{p^s} \sum_{\bb \in \Lambda_{p}}  |\dot{\cA_7}(\bb/p,M,-\b1)|^2 =
  \sum_{  \bgamma_1, \bgamma_2  \in  G_7 }          \frac{ \eta_M(\bgamma_1)\eta_M(\bgamma_2)
      \delta_{p}(  \bgamma_1 -   \bgamma_2  )}
     {\Nm(\bgamma_1)\Nm(\bgamma_2)   } 
\end{equation}
\begin{equation} \label{5.37a}
 = \sum_{\bb \in \Lambda_{p}}  \Big| \sum_{ \bgamma \in G_7, \;  \bgamma \equiv \bb \mod p}   \;\;       \frac{ \eta_M(\bgamma)}
      {\Nm(\bgamma)  } \Big|^2 . 
\end{equation}
Bearing in mind that $\eta_M(\bgamma)=0$ for $|\Nm(\bgamma)| \geq M $ (see (\ref{2.35a})), we get from (\ref{5.35}) that 
\begin{equation} \nonumber
  \varrho=  \sum_{\bb \in \Lambda_{p}}  \Big|  \sum_{\varsigma  =-1,1}  \sum_{\bk \in \cY_N}  \;\; 
       \sum_{ \substack{ \bgamma \in \Gamma^\bot \cap \cF, \;  \varsigma \bepsilon^{\bk}\bgamma  \in  \Gamma^\bot  \\
        \varsigma \bepsilon^{\bk}\bgamma \equiv \bb \mod p}}     \;\;       \frac{ \eta_M( \varsigma \bepsilon^{\bk}\bgamma)}
      {\Nm(\varsigma \bepsilon^{\bk}\bgamma)  } \Big|^2 . 
\end{equation}
We consider only $\bb = p_1\bb_0 \in \Lambda_p$, where $\bb_0 \in \Lambda_{p_2p_3}$ and  $p=p_1p_2p_3$. By  (\ref{2.0}), we obtain $\Gamma_{p_1 \fO} \subseteq  \Gamma^\bot  \subseteq \Gamma_{ \fO}$ and
 $\Gamma_{p_1\fO} = \{ \bgamma \in \Gamma^{\bot} | \gamma \equiv \bs \mod p_1\}$.
Hence,  we can take $\Gamma_{ p_1\fO}$ instead of $\Gamma^\bot$. 
We see $\varsigma \bepsilon^{\bk}\bgamma  \in  \Gamma_\fO$   for all $\bgamma \in \Gamma_\fO$, $\bk \in \ZZ^{s-1}$ and $\varsigma  \in \{-1,1\}$. 
Thus
\begin{equation} \nonumber
  \varrho \geq   \sum_{\bb \in \Lambda_{p_2p_3}}  \Big|\sum_{  \substack{ \varsigma  =-1,1 \\ \bk \in \cY_N}}  \;\; 
       \sum_{ \substack{ \bgamma \in \Gamma_\fO \cap \cF\\ 
        \varsigma \bepsilon^{\bk}\bgamma \equiv \bb \mod p_2p_3}}     \;\;       \frac{ \eta_M(p_1\varsigma \bepsilon^{\bk} \bgamma)}
      {\Nm(p_1\varsigma \bepsilon^{\bk}  \bgamma ) } \Big|^2 . 
\end{equation}
By  Lemma 4,  $(p_2,p_3) =1$. Hence, there exists $w_2, w_3 \in \ZZ$ such that \\ $p_2w_2 \equiv 1 \mod p_3$ and    $p_3w_3 \equiv 1 \mod p_2$.
It is easy to verify that if $\dot{\bb}_2, \ddot{\bb}_2  \in \Lambda_{p_2}$ (see (\ref{4.2a})), $\dot{\bb}_3, \ddot{\bb}_3  \in \Lambda_{p_3}$, and
$ (\dot{\bb}_2, \dot{\bb}_3) \neq   (\ddot{\bb}_2, \ddot{\bb}_3) $, then
\begin{equation} \nonumber
     \dot{\bb}_2 p_3w_3 +  \dot{\bb}_3 p_2w_2 \not\equiv   \ddot{\bb}_2 p_3w_3 +  \ddot{\bb}_3 p_2w_2 \mod  p_2p_3.
\end{equation}
Therefore
\begin{equation} \nonumber
   \Lambda_{p_2p_3} = \{ \bb \in \Lambda_{p_2p_3}    \; | \;   
   \;\;  \exists \;    \bb_2  \in \Lambda_{p_2}, \;    \bb_3  \in \Lambda_{p_3} \;\; 
        {\rm with}   \;\;   \bb \equiv  \bb_2 p_3w_3 +  \bb_3 p_2w_2 \mod  p_2p_3
    \}.
\end{equation}
Thus
\begin{equation} \nonumber
  \varrho \geq   \sum_{\bb_2 \in \Lambda_{p_2}}  \sum_{\bb_3 \in \Lambda_{p_3}}  \Big|\sum_{  \substack{ \varsigma  =-1,1 \\ \bk \in \cY_N}} \; 
       \sum_{ \substack{ \bgamma \in \Gamma_\fO \cap \cF\\ 
        \varsigma \bepsilon^{\bk}\bgamma  \equiv  \bb_2 p_3w_3 +  \bb_3 p_2w_2  \mod p_2p_3}}     \;\;       \frac{ \eta_M(p_1 \bgamma)}
      {\Nm(p_1\varsigma  \bgamma ) } \Big|^2  
\end{equation}
\begin{equation} \nonumber
 \geq   \sum_{\bb_2 \in \Lambda_{p_2}}  \sum_{\bb_3 \in \Lambda_{p_3}}  \Big|  \ddot{\chi}_{p_3}( \bb_3)
   \sum_{  \substack{ \varsigma  =-1,1 \\ \bk \in \cY_N}}  \; 
       \sum_{ \substack{ \bgamma \in \Gamma_\fO \cap \cF\\ 
        \varsigma \bepsilon^{\bk}\bgamma  \equiv  \bb_2 p_3w_3 +  \bb_3 p_2w_2  \mod p_2p_3}}     \;\;       \frac{ \eta_M(p_1 \bgamma)}
      {\Nm(p_1\varsigma  \bgamma ) } \Big|^2 
\end{equation}
\begin{equation} \nonumber
=   \sum_{\bb_2 \in \Lambda_{p_2}}  \sum_{\bb_3 \in \Lambda_{p_3}}  \Big|  \sum_{  \substack{ \varsigma  =-1,1 \\ \bk \in \cY_N}}  \; 
       \sum_{ \substack{ \bgamma \in \Gamma_\fO \cap \cF\\ 
        \varsigma \bepsilon^{\bk}\bgamma  \equiv  \bb_2 p_3w_3 +  \bb_3 p_2w_2  \mod p_2p_3}}     \;\;      
         \frac{ \ddot{\chi}_{p_3}(  \varsigma \bepsilon^{\bk}\bgamma)  \eta_M(p_1 \bgamma)}
      {\Nm(p_1\varsigma  \bgamma ) } \Big|^2 . 
\end{equation}
Using the Cauchy--Schwartz inequality, we have
\begin{equation} \nonumber
p_3^s  \varrho   \geq     \sum_{\bb_2 \in \Lambda_{p_2}}  \Big|  \sum_{\bb_3 \in \Lambda_{p_3}}  
  \sum_{  \substack{ \varsigma  =-1,1 \\ \bk \in \cY_N}}  \; 
       \sum_{ \substack{ \bgamma \in \Gamma_\fO \cap \cF\\ 
        \varsigma \bepsilon^{\bk}\bgamma  \equiv  \bb_2 p_3w_3 +  \bb_3 p_2w_2  \mod p_2p_3}}     \;\;      
         \frac{  \ddot{\chi}_{p_3}(  \varsigma \bepsilon^{\bk}\bgamma) \eta_M(p_1 \bgamma)}
      { p_1^{s} \Nm( \varsigma  \bgamma ) } \Big|^2 . 
\end{equation}
We see that  $ \varsigma \bepsilon^{\bk}\bgamma  \equiv  \bb_2 p_3w_3 \equiv \bb_2  \mod p_2$  if and only if there exists $\bb_3 \in \Lambda_{p_3}$ such that
$ \varsigma \bepsilon^{\bk}\bgamma  \equiv  \bb_2 p_3w_3 +  \bb_3 p_2w_2  \mod p_2p_3$.
Hence 
\begin{equation} \label{5.39}
 p_1^{2s}p_3^s \varrho   \geq     \sum_{\bb_2 \in \Lambda_{p_2}}  \Big|    \sum_{  \substack{ \varsigma  =-1,1 \\ \bk \in \cY_N}}  \; 
       \sum_{ \substack{ \bgamma \in \Gamma_\fO \cap \cF\\ 
        \varsigma \bepsilon^{\bk}\bgamma  \equiv  \bb_2  \mod p_2}}     \;\;      
         \frac{   \ddot{\chi}_{p_3}( \varsigma \bepsilon^{\bk}\bgamma) \eta_M(p_1 \bgamma)}
      {\Nm(\varsigma  \bgamma ) } \Big|^2 . 
\end{equation}
By (\ref{3.5}), we get $\Gamma_{i_0} = \varsigma \bepsilon^{\bk}   \Gamma_{i_0} $   for all $\bk \in \ZZ^{s-1}$,
 $\varsigma  \in \{-1,1\}$, and there exists $\Phi_{i_0} \subseteq \Lambda_{p_2}$ with
\begin{equation} \nonumber
    \Gamma_{i_0} = \bigcup_{\bb \in \Phi_{i_0}} (p_2 \Gamma_{\fO} + \bb), \quad  {\rm where} \quad 
    (p_2 \Gamma_{\fO} + \bb_1) \cap (p_2 \Gamma_{\fO} + \bb_2) = \emptyset, \; {\rm for} \; \; \bb_1 \neq \bb_2.
\end{equation} 
We consider in (\ref{5.39})   only $\bb_2 \in \Phi_{i_0}$. 
Applying the Cauchy--Schwartz inequality, we obtain
\begin{equation} \nonumber
  p_1^{2s} p_2^s p_3^s  \varrho   \geq    \Big|     \sum_{\bb_2 \in \Phi_{i_0}} \sum_{  \substack{ \varsigma  =-1,1 \\ \bk \in \cY_N}}  \; 
       \sum_{ \substack{ \bgamma \in \Gamma_\fO \cap \cF\\ 
        \varsigma \bepsilon^{\bk}\bgamma  \equiv  \bb_2  \mod p_2}}     \;\;      
         \frac{   \ddot{\chi}_{p_3}( \varsigma \bepsilon^{\bk}\bgamma) \eta_M(p_1 \bgamma)}
      {\Nm(\varsigma  \bgamma ) } \Big|^2  
\end{equation}
\begin{equation} \nonumber
=   \Big|    \sum_{  \substack{ \varsigma  =-1,1 \\ \bk \in \cY_N}}  \; 
       \sum_{  \bgamma \in \Gamma_{i_0} \cap \cF }     \;\;      
         \frac{   \ddot{\chi}_{p_3}( \varsigma \bepsilon^{\bk}\bgamma) \eta_M(p_1 \bgamma)}
      {\Nm(\varsigma  \bgamma ) } \Big|^2 . 
\end{equation}
Using Lemma 4, we get
\begin{equation} \nonumber
    \ddot{\chi}_{p_3}( \varsigma \bepsilon^{\bk}  \bgamma)  \frac{ |\Nm( \bgamma)|  }
      {\Nm(\varsigma  \bgamma)  }  =     \ddot{\chi}_{p_3}( \varsigma \bepsilon^{\bk}  \bgamma) \frac{ \Nm(\varsigma \bepsilon^{\bk}  \bgamma)  }
      {|\Nm(\varsigma \bepsilon^{\bk}  \bgamma)|  } 
\end{equation}
\begin{equation} \nonumber
 = \dot{\chi}_{p_3}( (\varsigma \bepsilon^{\bk}  \bgamma)) =
      \dot{\chi}_{p_3}( ( \bgamma)) =  \ddot{\chi}_{p_3}(  \bgamma) \frac{ |\Nm( \bgamma)|  }
      {  \Nm (\bgamma)  }.
\end{equation}
Hence
\begin{equation} \nonumber
   p_1^{2s}  p_2^s p_3^s  \varrho   \geq     \Big| \sum_{  \substack{ \varsigma  =-1,1 \\ \bk \in \cY_N}}   \;
       \sum_{ \bgamma \in \Gamma_{i_0} \cap \cF      }     \;\;       \frac{ \ddot{\chi}_{p_3}(  \bgamma) \eta_M(p_1\bgamma)}
      {\Nm(\bgamma) } \Big|^2 .
\end{equation}
Bearing in mind that $\eta_M(p_1\bgamma) = \eta_{M/p_1^s}(\bgamma)$ (see (\ref{2.35a})), we obtain
\begin{equation} \nonumber
    p_1^{2s}  p_2^s p_3^s \varrho \geq 4 \#\cY_N^2  \Big|        \sum_{ \bgamma \in \Gamma_{i_0} \cap \cF    }     \;\;       \frac{ \ddot{\chi}_{p_3}(  \bgamma) \eta_{M/p_1^s}(\bgamma)}
      {|\Nm(\bgamma)|  } \Big|^2 .
\end{equation}
Applying   Lemma 2, we have from  (\ref{5.36}) that $\#\cY_N \geq 0.5 c_5 (n/s)^{s-1}$ for $N \geq \dot{N}_0$ with some $\dot{N}_0>1$, and $n=s^{-1} \log_2 N$.
By Lemma 7 and  (\ref{5.37a}), we obtain 
\begin{equation} \nonumber
   \sup_{\bb \in \Lambda_{p}}  | \dot{\cA_7}(\bb/p,M,-\b1) |  \geq  2^{-s} \varrho^{1/2} \geq  c_{7} (2p_1^2 p_2p_3 )^{-s}  \#\cY_N  \geq  0.5c_{5} c_{7} (2p_1^2 p_2p_3 s)^{-s}   n^{s-1} ,
\end{equation}
with $M=[\sqrt{n}]=[\sqrt{\log_2 N}] \geq M_2 + \log_2  \dot{N}_0 $.
 Using Lemma 10, we get  the assertion of Lemma 12. \qed  \\
\\

{\bf 2.6. Auxiliary lemmas.} 
  We need the following notations  and results from [Skr]:\\

{\bf Lemma C.} [Skr, Lemma 3.2] {\it Let $\dot{\Gamma}  \subset \RR^s$ be an admissible lattice. Then
\begin{equation}\nonumber
        \sup_{\bx \in \RR^s}  \sum_{\bgamma \in \dot{\Gamma}} \prod_{1 \leq i \leq s} (1+|\gamma_i -x_i|)^{-2s}\leq H_{\dot{\Gamma}}
\end{equation} 
  where the constant $H_{\dot{\Gamma}}$  depends upon the lattice  $\dot{\Gamma}$ only by means of the invariants $\det \dot{\Gamma}$ and $\Nm \; \dot{\Gamma}$ }.\\

 Let $f(t), \; t\in \RR$, be a function of the class $C^{\infty}$; moreover let $f(t)$ and all derivatives $f^{(k)}$ belong to $L^1(\RR)$. 
We consider the following integrals for $\dot{\tau} >0$:
\begin{equation}\label{6.01}
  I(\dot{\tau}, \xi) = \int_{-\infty}^{\infty}{\frac{\eta(t)\widehat{\omega}(\dot{\tau} t) e(-\xi t)}{t}  dt}, \;   J_{f}(\dot{\tau}, \xi) = \int_{-\infty}^{\infty} f(t)\widehat{\omega}(\dot{\tau} t) e(-\xi t) dt.
\end{equation}

{\bf Lemma D.} [Skr, Lemma 4.2] {\it For all $\alpha>0$  and $\beta>0$, there exists a constant $ \breve{c}_{(\alpha,\beta)} >0$  such that
\begin{equation}\nonumber
 \max(|I(\dot{\tau}, \xi) |, |J_{f}(\dot{\tau}, \xi)|) <  \breve{c}_{(\alpha,\beta)} (1+\dot{\tau})^{-\alpha}(1+|\xi|)^{-\beta}.
\end{equation} 
}

 Let $m(t), \; t\in \RR$, be an even non negative function of the class $C^{\infty}$; moreover $m(t)=0$ for $|t| \leq 1$,  $m(t)=0$ 
 for $|t| \geq 4$,   and 
\begin{equation}\label{6.01a}
     \sum_{q= -\infty}^{+\infty} m(2^{-q}t) =1.
\end{equation} 
Examples of such functions see e.g. [Skr, ref. 5.16]. 
Let   $\dot{\bp}=(\dot{p}_1,...,\dot{p}_s)$, $\dot{p}_i>0, \; i=1,...,s$,  $a>0$, $x_0=\gamma_0=1$, 
\begin{equation}\label{6.02}
     \widehat{W}_{a,i}(\dot{\bp},\bx)  =  \frac{ \widehat{\omega} (\dot{p}_1 x_1) \eta(ax_1 )}{x_1}   \prod_{j=2}^s \frac{  
              \widehat{\omega} (\dot{p}_j x_j) m ( x_j)  }{x_j}  
    \frac{1}{ x_i} \quad {\rm for} \quad  \Nm \; \bx \neq 0,  
\end{equation}  
 and  $ \widehat{W}_{a,i}(\dot{\bp},\bx)=0$ for   $  \Nm ( \bx) =0  ,\;\;i=0,1,...,s$.    Let
\begin{equation}\label{6.03}
     \breve{W}_{a,i}(\dot{\Gamma},\dot{\bp},\bx)  =   \sum_{\bgamma\in \dot{\Gamma}^\bot \setminus {\bs} }    \widehat{W}_{a,i}(\dot{\bp},\bgamma)  e(\langle\bgamma,\bx\rangle).
\end{equation}        
By (\ref{2.10}) and (\ref{2.12}), we see that the series (\ref{6.03}) converge absolutely, and $ \widehat{W}_{a,i}(\dot{\bp},\bx)$ belongs to the class $C^{\infty}$. Therefore, we can use Poisson's summation formula (\ref{2.32}):
\begin{equation}\label{6.04}
    \breve{W}_{a,i}(\dot{\Gamma},\dot{\bp},\bx)  =  \det \dot{\Gamma}  \sum_{\bgamma\in \dot{\Gamma} }    W_{a,i}(\dot{\bp},\gamma - \bx) ,  
\end{equation} 
where $\widehat{W}_{a,i}(\dot{\bp},\bx)$ and $ W_{a,i}(\dot{\bp},\bx) $ are related by the Fourier transform. 
Using (\ref{6.02}),  we derive
\begin{equation}\nonumber
    W_{a,i}(\dot{\bp},\bx)  =   \prod_{j \in    \{1,...,s\}\setminus \{i\}}   w^{(1)}_1 (\dot{p}_j,x_j)  \prod_{j \in    \{1,...,s\} \cap \{i\}} 
   w^{(2)}_j (\dot{p}_j,x_j),
\end{equation} 
where co-factors can be described as follows (see also [Skr,  ref. 6.14-6.17]):

If $j=1$ and $i \neq 1$, then 
\begin{equation}\label{6.06}
 w^{(1)}_1(\tau, \xi) = \int_{-\infty}^{\infty}{\frac{1}{t} \eta(at)\widehat{\omega}(\tau t) e(-\xi t) dt} =I(a^{-1}\tau, a^{-1}\xi). 
\end{equation} 
Note that here we used formula (\ref{6.01}). If $j=1$ and $i = 1$, then 
\begin{equation}\nonumber
  w^{(2)}_1(\tau, \xi) = \int_{-\infty}^{\infty}{\frac{1}{t^2} \eta(at)\widehat{\omega}(\tau t) e(-\xi t) dt} =aJ_{f_1}(a^{-1}\tau, a^{-1}\xi). 
\end{equation}
Note that here we used formula (\ref{6.01})   with $f_1(t) =\eta(t)/t^2. $ If $j \geq 2$, then 
\begin{equation} \label{6.06a}
  w^{(l)}_j(\tau, \xi) = \int_{-\infty}^{\infty}{\frac{1}{t^l} m(t)\widehat{\omega}(\tau t) e(-\xi t) dt} =J_{f_2}(\tau, \xi). 
\end{equation}
Here we used formula  (\ref{6.01}) with $f_2(t) =m(t)/t^l, \quad j=2,...,s, \;l=1,2$.

Applying Lemma D, we obtain for $0<a \leq 1$ that
\begin{equation}\label{6.07}
|w_1^{(l)}(\tau, \xi) |  < \breve{c}_{(2s,2s)} (1+ a^{-1}|\xi|)^{-2s}, \; {\rm and }  \; |w_j^{(l)}(\tau, \xi)| <   \breve{c}_{(2s,2s)} (1+ |\xi|)^{-2s} , 
\end{equation} 
with $j=2,...,s, $ and $l=1,2$.
Now, using  (\ref{6.04}) and  Lemma C, we get (see also [Skr,  ref. 6.18, 6.19, 3.7, 3.10, 3,13]):\\

{\bf Lemma E.} {\it Let $\dot{\Gamma}  \subset \RR^s$ be an admissible lattice, and $0<a \leq 1$ . Then}
\begin{equation}\nonumber
        \sup_{\bx \in \RR^s} | \breve{W}_{a,i}(\dot{\Gamma},\dot{\bp},\bx)| \leq   \breve{c}_{(2s,2s)} \det \dot{\Gamma}  H_{\dot{\Gamma}}.
\end{equation} \\
\\

{\bf 2.7.  Dyadic decomposition of $\cB(\bb/p,M)$.} 
Using the  definition of the function $m(x)$ (see (\ref{6.01a})), we set
\begin{equation}  \label{6.1}
 \MM(\bx) = \prod_{j=2}^s m(x_j).
\end{equation}  
Let  $2^{\bq}=(2^{q_1},...,2^{q_s})$, and
\begin{equation}   \label{6.1a}
  \psi_{\bq}(\bgamma) =     \MM (2^{-\bq} \cdot \bgamma)\widehat{\Omega}(\tau  \bgamma)/\Nm(\bgamma),
\end{equation}
\begin{equation}\nonumber
  \cB_{\bq}(M)= \cB_{\bq}(\bb/p,M) =  \sum_{\bgamma\in \Gamma^\bot \setminus {\bs} }  \prod_{i=1}^s \sin(\pi \theta_i N_i \gamma_i)
    (  1-\eta_M(\bgamma)) 
           \psi_{\bq}(\bgamma)  )e(\langle\bgamma,\bb/p   \rangle  +\dot{x}),
\end{equation}
with $ \dot{x} = \sum_{1 \leq i  \leq s}  \theta_i N_i  \gamma_i/2  $. 

By (\ref{2.35}) and (\ref{6.01a}),  we have
\begin{equation}\label{6.2}
  \cB(\bb/p,M) =  \sum_{Q \in L} \cB_{\bq}(M),
\end{equation} 
with   $   \cL=\{\bq=(q_1,...,q_s) \in    \ZZ^s \; | \;  q_1+ \cdots +q_s =0\}$.

Let
\begin{equation}\label{6.8}
   \widetilde{\cB}_{\bq}(M) = \sum_{\bgamma\in \Gamma^\bot \setminus {\bs} } 
      \prod_{i=1}^s \sin(\pi \theta_i N_i \gamma_i)
      \eta(\gamma_1 2^{-q_1 }/M)
           \psi_{\bq}( \bgamma) e(\langle\bgamma,\bb/p  \rangle  +\dot{x}),
\end{equation}
and
\begin{equation}\nonumber
 \cC_{\bq}(M) =   \sum_{\bgamma\in \Gamma^\bot \setminus {\bs} }   \prod_{i=1}^s \sin(\pi \theta_i N_i \gamma_i)
       (  1-\eta_M(\bgamma)) 
         (1-   \eta(\gamma_1 2^{-q_1 }/M))   \psi_{\bq}( \bgamma)e(\langle\bgamma,\bb/p   \rangle  +\dot{x}).
 \end{equation}
According to (\ref{2.35a}), we get   $\eta_M(\bgamma) =1- \eta(2|\Nm(\bgamma)|/M)$,  $\eta(x)=0$ for $|x| \leq 1$, $\eta(x) = \eta(-x)$  and $\eta(x)=1$ for $|x| \geq 2$. 
Let $\eta(\gamma_1 2^{-q_1 }/M) m(\gamma_2 2^{-q_2 }) \cdots m(\gamma_s 2^{-q_s }) \neq 0$, then
$|\Nm(\bgamma)| \geq M $  (see (\ref{6.01a})), and
\begin{equation}\nonumber
                   (  1-\eta_M(\bgamma)) \eta(\gamma_1 2^{-q_1 }/M) = \eta(2|\Nm(\bgamma)|/M) \eta(\gamma_1 2^{-q_1 }/M)=\eta(\gamma_1 2^{-q_1 }/M).
\end{equation}
Hence
\begin{equation}\label{6.10}
               \cB_{\bq}(M) =   \widetilde{\cB}_{\bq}(M) +   \cC_{\bq}(M). 
\end{equation}
Let $n=s^{-1} \log_2 N, \; \tau =N^{-2}$ and
\begin{equation}\label{6.4}
   \cG_1=\{ \bq \in \cL \; |  \; \max_{i=1,...,s} q_i \geq -\log_2 \tau  + \log_2 n\},
\end{equation}
\begin{equation}\nonumber
   \cG_2=\{ \bq \in \cL \setminus \cG_1  \; |  \; \min_{i=2,...,s} q_i \leq - n-1/2\log_2 n \},
\end{equation}
\begin{equation}\nonumber
   \cG_3=\{ \bq \in \cL \; |  \; - n-1/2\log_2 n < \min_{i=2,...,s} q_i, \;\;  \max_{i=1,...,s} q_i < -\log_2 \tau  +\log_2 n  \},
\end{equation}
\begin{equation}\nonumber
   \cG_4=\{ \bq \in \cG_3 \; |  \;  q_1 \geq  -n  +s \log_2 n  \},
\end{equation}
\begin{equation}\nonumber
   \cG_5=\{ \bq \in \cG_3 \; |  \;   -n  -s \log_2 n  \leq   q_1 < -n  +s \log_2 n  \},
\end{equation}
\begin{equation}\nonumber
   \cG_6=\{ \bq \in \cG_3 \; |  \;   q_1 <   -n  -s \log_2 n   \}.
\end{equation}
We see
\begin{equation}\label{6.5}
              \cL   = \cG_1 \cup \cG_2 \cup \cG_3, \quad   \cG_3 = \cG_4 \cup \cG_5 \cup \cG_6 \quad {\rm and} \quad \cG_i \cap \cG_j =\emptyset, \; {\rm for}\; i\neq j
\end{equation}
and $i,j \in [1,3]$ or $i,j \in [4,6]$.   Let
\begin{equation}\label{6.6}
    \cB_i(M) = \sum_{{\bq} \in \cG_i}   \cB_{\bq}(M) .
\end{equation}  
By (\ref{6.2}), we obtain
\begin{equation}\label{6.3}
    \cB(\bb/p,M) =  \cB_1(M) + \cB_2(M) +\cB_3(M)  .
\end{equation} 
Let
\begin{equation}\label{6.11}
    \widetilde{\cB}_3(M)  = \sum_{{\bq} \in \cG_3}    \widetilde{\cB}_{\bq}(M), \quad   \widetilde{\cC}_3(M)  = \sum_{{\bq} \in \cG_3}   \cC_{\bq}(M).
\end{equation}
Applying (\ref{6.10}) and (\ref{6.6}), we get
\begin{equation}\label{6.3a}
           \cB_3(M) =   \widetilde{\cB}_3(M) +   \widetilde{\cC}_3(M).
\end{equation}
By (\ref{2.12}), we obtain the absolute convergence of the following series
\begin{equation}\nonumber
   \sum_{\bgamma\in \Gamma^\bot \setminus {\bs} } 
      | \widehat{\Omega} (\tau  \bgamma)/\Nm(\bgamma)|.
\end{equation}
Hence, the series (\ref{6.8}), (\ref{6.6}) and (\ref{6.11}) converges absolutely.

Let 
\begin{equation}\label{6.15}
    \breve{\cB}_{\bq}(M, \bvarsigma) =     \sum_{\bgamma\in \Gamma^\bot \setminus {\bs}}
           \eta(\gamma_1 2^{-q_1 }/M)\psi_{\bq}(\bgamma)e(\langle \bgamma, \bb/p+\dot{\btheta}(\bvarsigma) \rangle)
\end{equation}
with $\dot{\btheta}(\bvarsigma) = (\dot{\theta}_1(\bvarsigma),...,\dot{\theta}_s(\bvarsigma))$ and $\dot{\theta}_i(\bvarsigma) = (1+\varsigma_i) \theta_iN_i/4 , \; i=1,...,s$.
By (\ref{6.8}), we have 
\begin{equation} \label{6.15a}
         \widetilde{\cB}_{\bq}(M) =    \sum_{\bvarsigma \in \{1,-1\}^s }  \varsigma_1 \cdots \varsigma_s 
   (2\sqrt{-1})^{-s}  \breve{\cB}_{\bq}(M, \bvarsigma) .
\end{equation}
Let  $\bvarsigma_2 = -\b1 =-(1,1,...,1)$,  $\bvarsigma_3 = \dot{\b1}=(1,-1,...,-1)$, and let
\begin{equation}\label{6.18}
    \widetilde{\cB}_{3,1}(M)  = \sum_{{\bq} \in \cG_3}  \sum_{ \substack{\bvarsigma \in \{1,-1\}^s \\  \bvarsigma  \neq \bvarsigma_2, \bvarsigma_3 
    }}  \varsigma_1 \cdots \varsigma_s 
   (2\sqrt{-1})^{-s}  \breve{\cB}_{\bq}(M, \bvarsigma) , 
\end{equation}
\begin{equation}\label{6.19}
    \widetilde{\cB}_{i,j}(M)  = (-1)^{s+j} (2\sqrt{-1})^{-s} \sum_{{\bq} \in \cG_i}   \breve{\cB}_{\bq}(M, \bvarsigma_j), \quad i=3,4,5,6, \;j=2,3. 
\end{equation}
Using  (\ref{6.11}) and (\ref{6.15a}), we derive
\begin{equation}\nonumber
    \widetilde{\cB}_{3}(M)  = \widetilde{\cB}_{3,1}(M) +  \widetilde{\cB}_{3,2}(M) + \widetilde{\cB}_{3,3}(M) .
\end{equation}
Bearing in mind (\ref{6.5}), we obtain
\begin{equation}\label{6.23}
    \widetilde{\cB}_{3}(M)  = \widetilde{\cB}_{3,1}(M) +  \sum_{i=4,5,6}  \sum_{j=2,3} \widetilde{\cB}_{i,j}(M) .
\end{equation}
Let   
\begin{equation}\label{6.26}
       \widetilde{\cB}_{6,j,k}(M)  = (-1)^{s+j} (2\sqrt{-1})^{-s}  \sum_{{\bq} \in \cG_6}   \breve{\cB}_{\bq}^{(k)}(M,\bvarsigma_j ),  \quad  j=2,3, \;\; k=1,2,
\end{equation}
where
\begin{equation}\nonumber
   \breve{\cB}_{\bq}^{(1)}(M,\bvarsigma ) =     \sum_{\bgamma\in \Gamma^\bot \setminus {\bs}} 
           \eta(\gamma_1 2^{-q_1 }/M)\psi_{\bq}(\bgamma)\eta(2^{n +\log_2 n }\gamma_1)e(\langle \bgamma, \bb/p+\dot{\btheta}(\bvarsigma) \rangle)
\end{equation}
and
\begin{equation}\nonumber
   \breve{\cB}_{\bq}^{(2)}(M, \bvarsigma  ) =     \sum_{\bgamma\in \Gamma^\bot \setminus {\bs}} 
           \eta(\gamma_1 2^{-q_1 }/M)\psi_{\bq}(\bgamma)(1-\eta(2^{n +\log_2 n }\gamma_1))e(\langle \bgamma, \bb/p+\dot{\btheta}(\bvarsigma) \rangle).
\end{equation}
From  (\ref{6.15}), (\ref{6.19}) and  (\ref{6.26})  , we get
\begin{equation}\nonumber
         \breve{\cB}_{\bq}(M, \bvarsigma  )          = \breve{\cB}_{\bq}^{(1)}(M, \bvarsigma  )      +   \breve{\cB}_{\bq}^{(2)}(M, \bvarsigma  )
           \quad {\rm and} \quad  \widetilde{\cB}_{6,j}(M) =\widetilde{\cB}_{6,j,1}(M) + \widetilde{\cB}_{6,j,2}(M).
\end{equation}
So, we proved the following lemma:\\

{\bf Lemma 13.} {\it With notations as above, we get from (\ref{6.3}), (\ref{6.3a}) and (\ref{6.23})
\begin{equation}\label{6.28}
                    \cB(\bb/p,M) =  \bar{\cB}(M) +  \widetilde{\cC}_3(M) ,
\end{equation} 
where
\begin{equation} \label{6.27a}
               \bar{\cB}(M) =  \cB_{1}(M) +  \cB_{2}(M)  + \widetilde{\cB}_{3}(M) 
\end{equation}
and}
\begin{equation}\label{6.27}
    \widetilde{\cB}_{3}(M)  = \widetilde{\cB}_{3,1}(M) +  \sum_{j=2,3} (  \widetilde{\cB}_{4,j}(M)+\widetilde{\cB}_{5,j}(M)+\widetilde{\cB}_{6,j,1}(M)+\widetilde{\cB}_{6,j,2}(M)  ) .
\end{equation}

{\bf 2.8. The  upper bound estimate for $\bE(\bar{\cB}(M))$.} \\

{\bf Lemma 14.} {\it With notations as above
\begin{equation}\nonumber
   \cB_1(M) =O(1).
\end{equation}  }
{\bf Proof.}  Let $\bq \in \cG_1$, and let $j=q_{i_0} = \max_{1 \leq i \leq s}   q_i$, $i_0 \in [1,...,s]$. By  (\ref{6.4}), we have  $j \geq -\log_2 \tau + \log_2 n$.
Using (\ref{6.1a}), we obtain
\begin{equation}\label{7.2}
  | \cB_{\bq}(M)|  \leq \sum_{\bgamma \in \Gamma^\bot \setminus {\bs}} \Big| \prod_{i=1}^s \sin(\pi \theta_i N_i \gamma_i) 
     \frac{ \MM(2^{-\bq} \cdot\bgamma)\widehat{\Omega}(  \tau  \bgamma) }{  \Nm(\bgamma)   } \Big|.
\end{equation} 
From (\ref{6.1}) and (\ref{6.01a}), we get
\begin{equation}\label{7.2a}
    | \cB_{\bq}(M)|  \leq \rho_1 + \rho_2, \quad  {\rm with}  \quad \rho_i = \sum_{\bgamma \in \cX_i}  \frac{ |\MM(\bgamma)\widehat{\Omega}(\tau 2^{\bq}  \cdot \bgamma)| }{ |\Nm(\bgamma)|   },
\end{equation} 
where
\begin{equation}\nonumber
    \cX_1 =     \{ \bgamma \in 2^{-\bq} \cdot \Gamma^\bot \setminus {\bs}  \; | \;  |\gamma_1| \leq 2^{4sj}, \; |\gamma_i| \in [1,4], \; i =2,...,s  \} ,  
\end{equation} 
and
\begin{equation}\nonumber
    \cX_2 =     \{ \bgamma \in 2^{-\bq} \cdot \Gamma^\bot \setminus {\bs}  \; | \;  |\gamma_1| > 2^{4sj}, \; |\gamma_i| \in [1,4], \; i =2,...,s  \} .  
\end{equation} 
  We consider the admissible lattice $2^{-\bq} \cdot \Gamma^\bot$, where $\Nm (\Gamma^\bot) \geq 1$. Using  
  Theorem~A, we obtain that there exists a constant $c_{9}= c_{9}(\Gamma^\bot)$ such that
\begin{equation} \label{7.5}
   \# \{ \bgamma \in 2^{-\bq} \cdot \Gamma^\bot \; | \; |\gamma_i| \leq 4, i=2,...,s, \;2^{4(s-1)} |\gamma_1| \in [k,2k]  \} \leq c_{9} k,
\end{equation}
where $k=1,2,....$.

Let $i_0=1$.  We see that $\tau 2^{q_1} =\tau 2^{j}  \geq 2^{\log_2 n} =n$. By (\ref{2.12}),  (\ref{7.2}) and (\ref{7.5}), we get
\begin{equation}\nonumber
    \cB_{\bq}(M) = O(  \sum_{k \geq 0}    \sum_{\substack{\bgamma \in 2^{-\bq} \cdot \Gamma^\bot \setminus {\bs}, \;  1 \leq |\gamma_i| \leq 4,\;  i \geq 2\\ 2^{4(s-1)} |\gamma_1| \in [2^k,2^{k+1}]    }} 
    \frac{ |\widehat{\omega}(\tau 2^{q_1} \gamma_1)| }{ |\Nm(\bgamma)|   }= O\Big(  \sum_{k \geq 0}
     (1+ \tau 2^{q_1+k})^{-2s}  \Big).
\end{equation} 
Hence
\begin{equation}\label{7.7b}
   \cB_{\bq}(M) = O((\tau 2^{j})^{-2s}). 
\end{equation} 
Let $i_0 \geq 2$.  Bearing in mind (\ref{2.12}) and (\ref{7.5}), we have
\begin{equation}\nonumber
  \rho_1 = O\Big(  \sum_{0 \leq k \leq 4s(j+1)}   \sum_{\substack{\bgamma \in 2^{-\bq} \cdot \Gamma^\bot \setminus {\bs}, \;  1 \leq |\gamma_i| \leq 4,\;  i \geq 2\\ 2^{4(s-1)} |\gamma_1| \in [2^k,2^{k+1}]    }} 
    \frac{ |\widehat{\omega}(\tau 2^{q_{i_0}} \gamma_{q_{i_0}})| }{ |\Nm(\bgamma)|   }\Big)=
     O\Big(  \sum_{0 \leq k \leq 4s(j+1)} (1+ \tau  2^{q_{i_0}})^{-2s}  \Big).
\end{equation} 
Hence
\begin{equation}\label{7.6}
    \rho_1 =    O(j(1+ \tau 2^j)^{-2s}) .  
\end{equation} 
Taking into account that $q_1 =-(q_2+ \cdots + q_s) \geq -(s-1)j$ and $\tau 2^{j}  \geq n$, we obtain
\begin{equation}\nonumber
   \rho_2 = O\Big( \sum_{ k \geq 4sj}   \sum_{\substack{\bgamma \in 2^{-\bq} \cdot \Gamma^\bot \setminus {\bs}, \;  1 \leq |\gamma_i| \leq 4,\;  i \geq 2\\ 2^{4(s-1)} |\gamma_1| \in [2^k,2^{k+1}]    }} 
    \frac{ |\widehat{\omega}(\tau 2^{q_{1}} \gamma_{q_{1}})\widehat{\omega}(\tau 2^{q_{i_0}} \gamma_{q_{i_0}})| }{ |\Nm(\bgamma)|   }\Big)
\end{equation} 
\begin{equation}\nonumber
          = O\Big(  \sum_{k \geq 4sj} (1+ \tau 2^{q_1+k})^{-2s} (1+ \tau  2^{q_{i_0}})^{-2s}   \Big)         =O\Big( (1+ \tau  2^{q_{i_0}})^{-2s}  \Big).
\end{equation} 
Therefore
\begin{equation}\label{7.7}
    \rho_2   = O((1+ \tau 2^j)^{-2s}).  
\end{equation} 
Thus
\begin{equation}\label{7.7a}
   \cB_{\bq}(M) =O(j( \tau 2^j)^{-2s}).  
\end{equation} 
From  (\ref{3.2}), we have
\begin{equation}\label{7.8}
 \sum_{\bq \in \ZZ^{s},\; q_1+...+q_s=0,\; \max_{i} q_i =j} 1 =O(j^{s-2}).
\end{equation}
By  (\ref{6.4}), (\ref{6.6}), (\ref{7.7a}) and (\ref{7.7b}), we get
\begin{equation}\nonumber
   \cB_1(M) = \sum_{{\bq} \in \cG_1}   \cB_{\bq}(M)  =O\Big( \sum_{j \geq -\log_2 \tau+ \log_2 n}  \;\; \sum_{\bq \in \cL, \max_{i} q_i =j} j (\tau 2^j)^{-2s}  \Big)
\end{equation}
\begin{equation}\nonumber
 =O\Big( \sum_{j \geq -\log_2 \tau + \log_2 n}   j^{s}(\tau 2^j)^{-2s} \Big) =O(n^s(n)^{-2s}) =O(1).
\end{equation}
Hence, Lemma 14 is proved. \qed  \\

{\bf Lemma 15.}{\it With notations as above }
\begin{equation}\nonumber
  | \cB_2(M)| + |\widetilde{\cB}_{6,2,2}(M) +\widetilde{\cB}_{6,3,2}(M)|=O(n^{s-3/2}).
\end{equation}  
{\bf Proof.} We consider $\cB_2(M)$  (see (\ref{6.1a}), (\ref{6.4}) and (\ref{6.6})). 
Let $\bq \in \cG_2$, and let $j=-q_{i_0} = \min_{2 \leq i \leq s}   q_i$, $i_0 \in [2,...,s]$. 
 We see $j \geq n+1/2\log_2 n$ and $|\sin(\pi N_{i_0} \gamma_{i_0})| \leq \pi N_{i_0} 2^{-j+2}$ for $m(2^{-q_{i_0}} \gamma_{i_0} ) \neq 0$.
By (\ref{7.2}) and (\ref{7.2a}), we obtain
\begin{equation}\nonumber
   \cB_{\bq}(M) =O( \rho_1 + \rho_2), \quad  {\rm with}  \quad \rho_i = \sum_{\bq \in \cX_i}  \frac{ |N^{1/s}  2^{-j}\MM(\bgamma)\widehat{\Omega}(\tau 2^{\bq} \cdot  \bgamma)| }{ |\Nm(\bgamma)|   }.
\end{equation} 
Similarly to  (\ref{7.6}), (\ref{7.7}), we get 
\begin{equation}\nonumber
  \rho_1 = O\Big(  \sum_{0 \leq k \leq 4s(j+1)}   \sum_{\substack{\bgamma \in 2^{-\bq} \cdot \Gamma^\bot \setminus {\bs}, \;  1 \leq |\gamma_i| \leq 4,\;  i \geq 2\\ 2^{4(s-1)} |\gamma_1| \in [2^k,2^{k+1}]    }} 
    \frac{ N^{1/s}  2^{-j}}  { |\Nm(\bgamma)|   }\Big)
\end{equation} 
\begin{equation}\nonumber
    =       O\Big(  \sum_{0 \leq k \leq 4s(j+1)} N^{1/s}  2^{-j} \Big)=  O(jN^{1/s}  2^{-j}) .  
\end{equation} 
We see
\begin{equation}\nonumber
   \rho_2 = O\Big( \sum_{ k \geq 4sj}   \sum_{\substack{\bgamma \in 2^{-\bq} \cdot \Gamma^\bot \setminus {\bs}, \;  1 \leq |\gamma_i| \leq 4,\;  i \geq 2\\ 2^{4(s-1)} |\gamma_1| \in [2^k,2^{k+1}]    }} 
    \frac{ N^{1/s}  2^{-j}|\widehat{\omega}(\tau 2^{q_{1}} \gamma_{q_{1}})| }{ |\Nm(\bgamma)|   }\Big).
\end{equation} 
We have  $\max_{1 \leq i \leq s} q_i \leq - \log_2 \tau +\log_2 n$ for $\bq \in \cG_2$. Hence $q_1 = -(q_2+...+q_s) \geq (s-1)(\log_2 \tau -\log_2 n)$ and $\tau 2^{q_1} \geq \tau^{s} n^{-s+1} = 2^{-2ns} n^{-s+1} > 2^{-2sj}$.
 Thus
\begin{equation}\nonumber
 \rho_2 =    O\big(  N^{1/s}  2^{-j} \sum_{k \geq 4sj}  (1+ \tau 2^{q_1+k})^{-2s}   \big) =O\big( N^{1/s}  2^{-j}   \sum_{k \geq 4sj}  2^{-2s(k-2sj)}   \big) =O(N^{1/s}  2^{-j}).
\end{equation} 
Bearing in mind  (\ref{7.8}), we derive
\begin{equation}\nonumber
   \cB_2(M) = \sum_{{\bq} \in \cG_2}   \cB_{\bq}(M)  =O\Big( \sum_{j \geq n+1/2 \log_2 n} \sum_{\bq \in \cL, \min_{2 \leq i \leq s} q_i =-j} j N^{1/s}  2^{-j}  \Big)
\end{equation}
\begin{equation}\nonumber
 =O\Big( \sum_{j \geq n+1/2\log_2 n}   j^{s-1}N^{1/s}  2^{-j} \Big) =O(n^{s-3/2}).
\end{equation}
Consider  $\rho:=  \breve{\cB}_{\bq}^{(2)}(M,\dot{\b1} )+  \breve{\cB}_{\bq}^{(2)}(M,-\b1 )$. 
 By (\ref{6.1a}) and (\ref{6.26}),  we have
\begin{equation}\nonumber
   \rho =    O\Big(  \sum_{\bgamma\in \Gamma^\bot \setminus {\bs}} | \sin(\pi \theta_1 N_1 \gamma_1)
           \eta(\gamma_1 2^{-q_1 }/M)M (2^{-\bq} \cdot \bgamma)\widehat{\Omega}(\tau  \bgamma)/\Nm(\bgamma) 
\end{equation}
\begin{equation}\nonumber
   \times   (1 - \eta(2^{n +\log_2 n }\gamma_1)) e(\langle \bgamma, \bb/p \rangle) |\Big)
\end{equation}
\begin{equation}\nonumber
  =O\Big(   \sum_{\bgamma\in 2^{-\bq} \dot\Gamma^\bot \setminus {\bs}} |\sin(\pi \theta_1 N_1 2^{q_1}\gamma_1) (1 - \eta(2^{q_1+n +\log_2 n }\gamma_1))
          \MM (\bgamma)/\Nm(\bgamma) | \Big).
\end{equation} 
Applying    (\ref{2.35a}), (\ref{6.1}) and (\ref{7.5}), we obtain
\begin{equation}\nonumber
   \rho =O\Big(   \sum_{\bgamma\in 2^{-\bq}\Gamma^\bot \setminus {\bs},\; |\gamma_1| \leq 2^{-q_1-  n -\log_2 n +4 }} | N_1 2^{q_1}\gamma_1 
          \MM ( \bgamma)/\Nm(\bgamma) | \Big) =O(1/n).
\end{equation}
\begin{equation}\nonumber
  =O\Big(   \sum_{ \substack{\bgamma \in 2^{-\bq} \cdot \Gamma^\bot \setminus {\bs}, \;  1 \leq |\gamma_i| \leq 4,\;  i \geq 2\\ |\gamma_1| \leq 2^{-q_1-  n -\log_2 n +4 }}}  N_1 2^{q_1}  \Big) 
  =O(N_1 2^{q_1}2^{-q_1-  n -\log_2 n +4 } )
  =O(1/n).
\end{equation}
We get from (\ref{6.4}) that 
\begin{equation}\label{7.17a}
                     \# \cG_3 =O(n^{s-1}).
\end{equation}
By (\ref{6.4}) and (\ref{6.26}),  we get $ \widetilde{\cB}_{6,2,2}(M) +\widetilde{\cB}_{6,3,2}(M)  =(n^{s-2})$.
 
Hence, Lemma 15 is proved. \qed  \\

{\bf Lemma 16.} {\it With notations as above}
\begin{equation}\nonumber
      | \bE(\widetilde{\cB}_{3,1}(M))|  +  | \bE(\widetilde{\cB}_{4,3}(M))|+ | \widetilde{\cB}_{5,2}(M)| +| \widetilde{\cB}_{5,3}(M)|=  O(n^{s-3/2}).
\end{equation}
{\bf Proof.}  By (\ref{6.1a}) and (\ref{6.15}), we have
%
\begin{equation}\label{7.27}
    \breve{\cB}_{\bq}(M, \bvarsigma) =     \sum_{\bgamma\in  2^{-\bq} \cdot \Gamma^\bot \setminus {\bs}}
           \eta(\gamma_1 /M)\psi_{\bq}(2^{\bq} \cdot\bgamma)e(\langle \bgamma,\bx\rangle)
\end{equation}
\begin{equation} \nonumber
     =    \sum_{\bgamma\in 2^{-\bq} \cdot \Gamma^\bot \setminus {\bs}}       
          \frac{ \widehat{\omega} (2^{q_1}\tau \gamma_1) \eta(\gamma_1 /M)}{\gamma_1}   \prod_{j=2}^s \frac{  
              \widehat{\omega} (2^{q_j}\tau \gamma_j) m ( \gamma_j)  }{\gamma_j}
     e(\langle\bgamma,\bx \rangle)   ,
\end{equation}
with $\bx= 2^{\bq} \cdot( \bb/p+\dot{\btheta}(\bvarsigma))  $ and $\dot{\theta}_i(\bvarsigma) = (1+\varsigma_i) \theta_iN_i/4 , \; i=1,...,s$.\\ 
Applying (\ref{6.04})  and Lemma E  with $\dot{\Gamma} = 2^{-\bq} \Gamma$, $i =0$, and $\dot{\bp} = \tau 2^{\bq}$, we get
\begin{equation}  \nonumber 
           \breve{\cB}_{\bq}(M,\bvarsigma ) = O(   1 ).
\end{equation}
Using  (\ref{6.4}), we  obtain   $\#  \cG_5 =O(n^{s-2} \log_2 n) $. \\
By (\ref{6.19}) , we get
\begin{equation}\label{7.27a}
    \widetilde{\cB}_{5,i}(M)   =  O\Big(  
       \sum_{\bq \in  \cG_5}  |\breve{\cB}_{\bq}(M,\bvarsigma)|    \Big)  = O(n^{s-2} \log_2 n),\quad i=2,3 .
\end{equation}
Consider $\bE(\widetilde{\cB}_{3,1}(M))$ and $\bE(\widetilde{\cB}_{4,3}(M))$. Let 
\begin{equation}\nonumber
        \bE_i(f) =\int_{0}^1 f(\btheta)d\theta_i.      
\end{equation} 
Let $\bvarsigma \neq -\b1 $. Then  there exists $i_0 =i_0(\bvarsigma) \in [1,s]$  with $\varsigma_{i_0} =1$.\\
 By  (\ref{5.20a})  and (\ref{7.27}),  we have
\begin{equation}\nonumber
       \bE_{i_0}(\widetilde{\cB}_{\bq}(M,\bvarsigma ))  =    \sum_{\bgamma\in  2^{-\bq} \cdot \Gamma^\bot \setminus {\bs}}   
           \frac{e(N_{i_0}  2^{q_{i_0}} \gamma_{i_0}/2)-1}{\pi \sqrt{-1} N_{i_0}  2^{q_{i_0}} \gamma_{i_0}}  \;  \frac{ \widehat{\omega} (2^{q_1} \tau \gamma_1) \eta(\gamma_1 /M)}{\gamma_1} 
\end{equation}
\begin{equation}\nonumber
   \times     \prod_{j=2}^s \frac{  
              \widehat{\omega} ( 2^{q_j} \tau \gamma_j) m ( \gamma_j)  }{\gamma_j}  e(\langle\bgamma,\bx\rangle),
\end{equation}
with some $\bx \in \RR^s$. 
Hence
\begin{equation}\nonumber
       \bE_{i_0}(\breve{\cB}_{\bq}(M,\bvarsigma )) = O\Big(   N_{i_0}^{-1} 2^{-q_{i_0}}   \sup_{\bx \in \RR^s} \Big|  \sum_{\bgamma\in 2^{-\bq} \cdot \Gamma^\bot \setminus {\bs}}  
       \widehat{\cB}_{\bq}(M, \bgamma, i_0)     
        e(\langle\bgamma,\bx \rangle)   \Big|    \Big),
\end{equation}
where
\begin{equation}\nonumber
     \widehat{\cB}_{\bq}(M, \bgamma,i_0)     =  \frac{ \widehat{\omega} (2^{q_1}\tau \gamma_1) \eta(\gamma_1 /M)}{\gamma_1}   \prod_{j=2}^s \frac{  
              \widehat{\omega} (2^{q_j}\tau \gamma_j) m ( \gamma_j)  }{\gamma_j}
     \frac{1}{\gamma_{i_0}}  .  
\end{equation}
Applying (\ref{6.04})  and Lemma E  with $\dot{\Gamma} = 2^{-\bq} \Gamma$,  and $\dot{\bp} = \tau 2^{\bq}$, we obtain
\begin{equation}\label{7.23}
       \bE(\breve{\cB}_{\bq}(M,\bvarsigma )) =     \bE(\bE_{i_0}(\breve{\cB}_{\bq}(M,\bvarsigma ))) = O(    N_{i_0}^{-1} 2^{-q_{i_0}}   ).
\end{equation}
By  (\ref{6.18}), we have $i_0(\bvarsigma) \geq 2$ and 
\begin{equation}\nonumber
       \bE(\widetilde{\cB}_{3,1}(M))    =  O\Big(  \sum_{ \substack{\bvarsigma \in \{1,-1\}^s  \\  \bvarsigma  \neq -\b1, \dot{\b1} }} 
       \sum_{\bq \in  \cG_3}    N_{i_0(\bvarsigma)}^{-1} 2^{-q_{i_0(\bvarsigma)}}  \Big).
\end{equation}
Using (\ref{6.4}), we get $\# \{\bq \in \cG_3\; |\; q_{i_0} =j   \}=O(n^{s-2} ) $ and $ j \geq - n - 1/2\log_2 n$. Hence
\begin{equation}\label{7.23a}
        \bE(\widetilde{\cB}_{3,1}(M))   =  O\Big( n^{s-2}  \sum_{  j \geq - n - 1/2\log_2 n}     N^{-1/s} 2^{-j} \Big) =O(n^{s-3/2} ).
\end{equation}
From (\ref{6.4}), we get $q_1 \geq  -n +s\log_2 n$ for $\bq\in  \cG_4$. Applying (\ref{6.19}),
 (\ref{7.17a}) and (\ref{7.23}) with $i_0(\bvarsigma) =1$, we obtain   
\begin{equation}\nonumber
        \bE(\widetilde{\cB}_{4,3}(M))   =O\Big( \sum_{\bq \in  \cG_4}  N_{1}^{-1} 2^{-q_{1}}\Big) = O\Big( n^{s-1}  \sum_{  q_1 \geq  -n +s\log_2 n}    N^{-1/s} 2^{-q_1} \Big) =O(1).
\end{equation}
By (\ref{7.27a}) and (\ref{7.23a}), Lemma 16 is proved. \qed \\

{\bf Lemma 17.} {\it With notations as above}
\begin{equation}\nonumber
     \widetilde{\cB}_{4,2}(M) =  O( n^{s-3/2} ).
\end{equation}
{\bf Proof.} By (\ref{7.27}), we have
\begin{equation} \nonumber
    \breve{\cB}_{\bq}(M,-\b1) =    \sum_{\bgamma\in 2^{-\bq} \cdot \Gamma^\bot \setminus {\bs}}       
          \frac{ \widehat{\omega} (2^{q_1}\tau \gamma_1) \eta(\gamma_1 /M)}{\gamma_1}   \prod_{j=2}^s \frac{  
              \widehat{\omega} (2^{q_j}\tau \gamma_j) m ( \gamma_j)  }{\gamma_j}
     e(\langle\bgamma,2^{\bq} \cdot \bb/p \rangle)   .
\end{equation}
From (\ref{6.06}), we derive that $I(d,v)=0$ for $v=0$. Hence $ w^{(1)}_1(\tau, 0)=0$. Now applying (\ref{6.04}) - (\ref{6.07}) with $\dot{\Gamma}^\bot= 2^{-\bq} \cdot \Gamma^\bot, i=0$ and $a=M^{-1}$,    we get
\begin{equation}\nonumber
      |\breve{\cB}_{\bq}(M,-\b1) |  \leq   \breve{c}_{(2s,2s)} \det \Gamma  \sum_{\bgamma\in 2^{\bq} \cdot
       \Gamma,
      \; \gamma_1 \neq (\bb/p)_1  } (1+M|\gamma_1 -2^{q_1}(\bb/p)_1 |)^{-2s}
\end{equation}
\begin{equation}\nonumber
    \times    \prod_{i=2}^s (1+|\gamma_i -2^{q_i}(\bb/p)_i|)^{-2s}.   
\end{equation}
Bearing in mind (\ref{2.0}), we  get $p_1 \Gamma_{\fO}  \subseteq \Gamma^\bot \subseteq \Gamma_{\fO}$.
 Taking into account that $p=p_1p_2p_3$ and $\bb \in \Gamma_{\fO}$, we obtain 
\begin{equation}\label{7.35}
      |\breve{\cB}_{\bq}(M,-\b1) |  \leq   \breve{c}_{(2s,2s)} \det \Gamma p^{2s^2}\sum_{\bgamma\in   p  2^{\bq} \cdot \Gamma \setminus \bs  } (1+M|\gamma_1  |)^{-2s}   \prod_{i=2}^s (1+|\gamma_i|)^{-2s}.   
\end{equation}
We have
\begin{equation}\label{7.36}
      |\breve{\cB}_{\bq}(M,-\b1) |  \leq  \breve{c}_{(2s,2s)} \det \Gamma  p^{2s^2}(a_1+a_2), 
\end{equation}
where
\begin{equation}\nonumber
  a_1 = \sum_{\bgamma \in   p  2^{\bq} \cdot \Gamma \setminus \bs , \max |\gamma_i|  \leq M^{1/s }} (1+M|\gamma_1  |)^{-2s}   \prod_{i=2}^s (1+|\gamma_i|)^{-2s},
\end{equation}
and
\begin{equation}\nonumber
  a_2 = \sum_{\bgamma \in   p  2^{\bq} \cdot\Gamma \setminus \bs , \max |\gamma_i| >M^{1/s}} (1+M|\gamma_1  |)^{-2s}   \prod_{i=2}^s (1+|\gamma_i|)^{-2s}.   
\end{equation}
We see that $|\gamma_1| \geq M^{-(s-1)/s}$ for $\max_{1 \leq i \leq s} |\gamma_i|  \leq M^{1/s}$.
Applying Theorem A, we have 
\begin{equation}\nonumber
  a_1 \leq  M^{-2} \sum_{\bgamma \in   p  2^{\bq} \cdot\Gamma \setminus \bs , \max |\gamma_i|  \leq M^{1/s }}   1  = O(M^{-1}),
\end{equation}
and
\begin{equation}\label{7.40}
  a_2 \leq  \sum_{j \geq M^{1/s }}  \sum_{  \substack{ \bgamma \in   p  2^{\bq} \cdot\Gamma \setminus \bs \\ \max |\gamma_i| \in [j,j+1) }} j^{-2s} = O\Big(\sum_{j \geq M^{1/s } } j^{-s}\Big) =O(M^{-(s-1)/s}).
\end{equation}
Taking into account that $\# \cG_3 =O(n^{s-1})$ (see (\ref{7.17a})), we get from (\ref{7.36})  and  (\ref{6.19})   that
\begin{equation}\nonumber
    \widetilde{\cB}_{4,2}(M)  =  O\Big(         \sum_{\bq \in  \cG_4}  \breve{\cB}_{\bq}(M,-\b1) \Big)= 
     O\Big(         \sum_{\bq \in  \cG_3}  M^{-1/2} \Big)  = O(M^{-1/2} n^{s-1} ).
\end{equation}
Hence, Lemma 17 is proved. \qed  \\

{\bf Lemma 18.} {\it With notations as above}
\begin{equation}\nonumber
   \widetilde{\cB}_{6,2,1}(M) + \widetilde{\cB}_{6,3,1}(M) =  O( n^{s-3/2}), \quad M=[\sqrt{n}].
\end{equation}
{\bf Proof.} Let  $M_1=2^{-q_1-n -\log_2 n }$. By (\ref{6.4}), we get $M_1 \geq n \geq 2M$ for $\bq \in \cG_6$ and $n \geq 4$. From (\ref{2.35a}), we  have $  \eta(\gamma_1/M) \eta(\gamma_1/M_1)= \eta(\gamma_1/M_1)$.
 Using (\ref{6.1a}),  (\ref{6.15}) and (\ref{6.26}), we derive similarly to (\ref{7.27}) that
\begin{equation}\nonumber
   \breve{\cB}_{\bq}^{(1)}(M,\bvarsigma_j ) =     \sum_{\bgamma\in 2^{\bq}\cdot \Gamma \setminus {\bs}}    \frac{ \widehat{\omega} (2^{q_1}\tau \gamma_1) \eta(\gamma_1 /M_1)}{\gamma_1}  
\end{equation}
\begin{equation}\nonumber
   \times   \prod_{i=2}^s  \frac{ \widehat{\omega} (2^{q_j}\tau \gamma_j) m ( \gamma_j)  }{\gamma_j}  e(\langle \bgamma, 2^{\bq}\cdot( \bb/p +(j-2)\theta_1N_1(1,0,...,0))\rangle)
\end{equation}
with $j=2,3$, $\bvarsigma_2=-\b1$ and $\bvarsigma_3=\dot{\b1}$. \\
By(\ref{6.06a}), we obtain that, 
 $J_{f_2}(\tau, v)=0$  with $f_2(t) =m(t)/t$ for $v=0$ . Hence $ w^{(1)}_2(\tau, 0)=0$. Now applying (\ref{6.04}) - (\ref{6.07}) with $\dot{\Gamma}^\bot= 2^{-\bq} \cdot \Gamma^\bot$, $i=0$ and
 $a=M_1^{-1}=2^{q_1+n +\log_2 n }$,     we get analogously to (\ref{7.35})   
\begin{equation}\nonumber
      |\breve{\cB}_{\bq}^{(1)}(M,\bvarsigma_j )|  \leq   \breve{c}_{(2s,2s)} \det \Gamma  p^{2s^2}\sum_{\bgamma\in   p  2^{\bq} \cdot\Gamma \setminus \bs  } (1+M_1|\gamma_1 -x(j) |)^{-2s}   \prod_{i=2}^s (1+|\gamma_i|)^{-2s},   
\end{equation}
 with $x(j)= (j-2)p\theta_12^{q_1}N_1$.  We have
\begin{equation}\label{7.46}
      |  \breve{\cB}_{\bq}^{(1)}(M,\bvarsigma_j) |  \leq  \breve{c}_{(2s,2s)} \det \Gamma p^{2s^2}(a_3+a_4), 
\end{equation}
where
\begin{equation}\nonumber
  a_3 = \sum_{\bgamma \in   p  2^{\bq} \cdot\Gamma\setminus \bs, \;\max |\gamma_i|  \leq M^{1/s }} (1+M_1|\gamma_1 -x(j)|)^{-2s}   \prod_{i=2}^s (1+|\gamma_i|)^{-2s},
\end{equation}
and
\begin{equation}\nonumber
  a_4 = \sum_{\bgamma \in   p  2^{\bq} \cdot\Gamma, \max |\gamma_i| >M^{1/s}} (1+M_1|\gamma_1 -x(j) |)^{-2s}   \prod_{i=2}^s (1+|\gamma_i|)^{-2s}.   
\end{equation}
We see that $|\gamma_1| \geq M^{-(s-1)/s}$ for $\max_{1 \leq i \leq s} |\gamma_i|  \leq M^{1/s}$.
Bearing in mind that \\$|x(j)| \leq c_3p n^{-s}$ for $\bq \in \cG_6$, we obtain $|\gamma_1|  \geq 2|x(j)|$
 for $M=[\sqrt{n}]$ and $N> 8psc_3$. 
Applying Theorem A, we get 
\begin{equation}\nonumber
  a_3 \leq 2^{2s} M_1^{-2s} M^{2(s-1)}\sum_{\bgamma \in   p  2^{\bq} \cdot\Gamma, \; \max |\gamma_i|  \leq M^{1/s }}   1  = O(M^{-1}).
\end{equation}
Similarly to (\ref{7.40}), we have
\begin{equation}\nonumber
  a_4 \leq  \sum_{j \geq M^{1/s }}  \sum_{  \substack{ \bgamma \in   p  2^{\bq} \cdot\Gamma \setminus \bs \\ \max |\gamma_i| \in [j,j+1) }} j^{-2s} = O\big(\sum_{j \geq M^{1/s } } j^{-s}\big) =O(M^{-(s-1)/s}).
\end{equation}
By (\ref{6.4}) and (\ref{7.17a}), we obtain
$\# \cG_6   \leq \# \cG_3 =O(n^{s-1})$.  We get from (\ref{6.26}) and (\ref{7.46}) that 
\begin{equation}\nonumber
     \widetilde{\cB}_{6,2,1}(M) + \widetilde{\cB}_{6,3,1}(M)   =  O\Big(         \sum_{\bq \in  \cG_6, \; j=2,3}     \breve{\cB}_{\bq}^{(1)}(M,\bvarsigma_j ) \Big)  = O(M^{-1/2} n^{s-1}).
\end{equation}
Hence, Lemma 18 is proved. \qed \\

Using (\ref{6.27}), (\ref{6.27a})  and Lemma 14 - Lemma 18, we obtain \\
{\bf  Corollary 1.} {\it With notations as above}
\begin{equation}\nonumber
       \bE(\bar{\cB}(M)) =  O( n^{s-5/4}), \quad  M=[\sqrt{n}].
\end{equation}  \\

{\bf 2.9. The  upper bound estimate for $\bE({ \widetilde{\cC}_3}(M))$ and Koksma--Hlawka inequality.}
Let
\begin{equation}\nonumber
   \cG_7=\{ \bq \in \cG_3 \; |  \;   -\log_2 \tau  -s\log_2 n  \leq   \max_{i=1,...,s} q_i < -\log_2 \tau  +\log_2 n  \}.
\end{equation}
\begin{equation}\label{9.1}
   \cG_8=\{ \bq \in   \cG_3 \setminus \cG_7 \; |  \;    q_1 < -n -1/2 \log_2 n  \},
\end{equation}
\begin{equation}  \nonumber
   \cG_9=\{ \bq \in   \cG_3 \setminus \cG_7 \; |  \;    q_1 \geq -n  -1/2 \log_2 n  \},
\end{equation}
and let
\begin{equation}\nonumber
     \widetilde{\cC}_i(M) = \sum_{{\bq} \in \cG_i}   \cC_{\bq}(M) , \quad i=7,8,9.
\end{equation}  
It is easy to see that 
\begin{equation}  \nonumber
  \cG_3 = \cG_7 \cup \cG_8 \cup \cG_9,  \quad {\rm and}\;   \quad \cG_i \cap \cG_j =\emptyset, \quad  {\rm for}\; i \neq j.
\end{equation} 
Hence
\begin{equation}\label{9.3a}
  \widetilde{\cC}_3(M) =   \widetilde{\cC}_7(M) +  \widetilde{\cC}_8(M) + \widetilde{\cC}_9(M)  .
\end{equation}
From (\ref{6.8}), we have similarly to (\ref{6.15}) that 
\begin{equation}\label{9.5}
         \cC_{\bq}(M) =    \sum_{\bvarsigma \in \{1,-1\}^s}  \varsigma_1 \cdots \varsigma_s 
  (2\sqrt{-1})^{-s}  \breve{\cC}_{\bq}(M, \bvarsigma) ,
\end{equation}
where
\begin{equation}\nonumber
    \breve{\cC}_{\bq}(M, \bvarsigma) =     \sum_{\bgamma\in \Gamma^\bot \setminus {\bs}}
          \psi_{\bq}(\bgamma) (  1-\eta_M(\bgamma))(1-   \eta(\gamma_1 2^{-q_1 }/M))  e(\langle \bgamma, \bb/p+\dot{\btheta}(\bvarsigma)\rangle),
\end{equation}
with $\dot{\theta}_i(\bvarsigma) = (1+\varsigma_i) \theta_iN_i/4 , \; i=1,...,s$.

By (\ref{9.5}) and (\ref{9.1}), we get  
\begin{equation}\label{9.8}
     \widetilde{\cC}_{9}(M)  =  \widetilde{\cC}_{10}(M) +  \widetilde{\cC}_{11}(M) ,
\end{equation}
where
\begin{equation}\label{9.6}
     \widetilde{\cC}_{10}(M)  = \sum_{{\bq} \in \cG_9}  \sum_{ \substack{\bvarsigma \in \{1,-1\}^s \\  \bvarsigma   \neq -\b1}}  \varsigma_1 \cdots \varsigma_s 
  (2\sqrt{-1})^{-s}  \breve{\cC}_{\bq}(M, \bvarsigma) , 
\end{equation}
and
\begin{equation}\label{9.7}
     \widetilde{\cC}_{11}(M)  = (-1)^s   (2\sqrt{-1})^{-s}  \sum_{{\bq} \in \cG_9}   \breve{\cC}_{\bq}(M, -\b1). 
\end{equation} \\

{\bf Lemma 19.} {\it With notations as above}
\begin{equation}  \nonumber
     \bE( \widetilde{\cC}_{i}(M))    =    O(n^{s-3/2}), \quad  i=7,8,10, \quad M =[\sqrt n].
\end{equation}  
{\bf Proof.}   Let $\bgamma \in 2^{-\bq} \cdot \Gamma^\bot \setminus {\bs}  $. By (\ref{2.35a}), (\ref{6.01a}) and  (\ref{6.1}), we have \\$ (  1-\eta_M(\bgamma))(1-   \eta(\gamma_1 /M)) \MM(\bgamma) \neq 0$ only if $ 2^{-2s+3} M \leq |\gamma_1| \leq   2M, \; |\gamma_i| \in [1,4]$, $i =2,...,s   $.
From(\ref{6.8}), we derive
\begin{equation}\label{9.10}
   \cC_{\bq}(M) =O\Big(  \sum_{\bgamma \in \cX}  \Big|\prod_{i=1}^s \sin(\pi \theta_i N_i 2^{q_i} \gamma_i) \frac{ \MM(\bgamma)\widehat{\Omega}( \tau 2^{\bq} \cdot  \bgamma) }{ \Nm(\bgamma)  } \Big| \Big)
\end{equation} 
where
\begin{equation}\nonumber
    \cX =     \{ \bgamma \in 2^{-\bq} \cdot \Gamma^\bot \setminus {\bs}  \; | \; 2^{-2s+3} M \leq |\gamma_1| \leq   2M, \; |\gamma_i| \in [1,4], \; i =2,...,s  \} . 
\end{equation}
Bearing in mind (\ref{7.5}), we get $ \cC_{\bq}(M) =O( 1 )$.

Using (\ref{3.2}), (\ref{6.4}) and (\ref{9.1}),   we obtain  $\#  \cG_7 =O(n^{s-2} \log_2 n) $.
Applying(\ref{9.1}), we get
\begin{equation}\label{9.11}
   \widetilde{\cC}_{7}(M)  = \sum_{{\bq} \in \cG_7}   \cC_{\bq}(M)  = O( n^{s-2} \log_2 n).
\end{equation} 
 Consider $ \widetilde{\cC}_{8}(M)$. Let $\bgamma \in \cX$. Then   $|\sin(\pi \theta_1 N_{1} 2^{q_1}\gamma_{1})| \leq  \pi M N_{1} 2^{1+q_1}$.
 
By (\ref{9.10}), we have
\begin{equation}\nonumber
   \cC_{\bq}(M) =O\Big( \sum_{\bgamma \in \cX}  \frac{ |MN^{1/s}   2^{q_1}\widehat{\Omega}(\tau 2^{\bq} \cdot  \bgamma)| }{ |\Nm(\bgamma)|   }\Big)
      = O( M N^{1/s}  2^{q_1}).
\end{equation} 
Using (\ref{3.2}) and (\ref{9.1}),   we derive  $\#  \{ \bq \in \cG_8 | q_1=d\} =O(n^{s-2} ) $. Hence
\begin{equation}\nonumber
    \widetilde{\cC}_8(M) = \sum_{{\bq} \in \cG_8}   \cC_{\bq}(M)  =O\Big( \sum_{j \geq n+ 0.5 \log_2 n} \;
    \sum_{{\bq} \in \cG_8, \; q_1 =-j} M N^{1/s}   2^{-j}  \Big)
\end{equation}
\begin{equation}\label{9.14}
 =O\Big(n^{s-2} M \sum_{j \geq n +  0.5\log_2 n}    2^{n-j} \Big) =O(n^{s-2}).
\end{equation}

 Consider $ \widetilde{\cC}_{10}(M)$. 
From(\ref{9.6}),  we get that there exists $i_0 =i_0(\bvarsigma) \in [1,s]$  with $\varsigma_{i_0} =1$.
 By  (\ref{5.20a}),  (\ref{6.1a})  and (\ref{9.5}),  we have
\begin{equation}\nonumber
       \bE_{i_0}(\breve{\cC}_{\bq}(M,\bvarsigma ))  =    \sum_{\bgamma\in \Gamma^\bot \setminus {\bs}}   \dot{\cC}_{\bq}(M, \bgamma)  
           \frac{e(N_i \gamma_{i_0}/2)-1}{\pi\sqrt{-1} N_{i_0} \gamma_{i_0}}   
        e(\langle\bgamma,\bx\rangle)  
\end{equation}
with some  $\bx  \in \RR^s$, where
\begin{equation}\nonumber
    \dot{\cC}_{\bq}(M, \bgamma) =   
            (  1-\eta_M(\bgamma))(1-   \eta(\gamma_1 2^{-q_1 }/M))   \widehat{\Omega} ( \tau \cdot \bgamma) \MM(2^{-\bq}  \bgamma)  / \Nm(\bgamma)  .
\end{equation}
Hence
\begin{equation}\nonumber
       \bE_{i_0}(\breve{\cC}_{\bq}(M,\bvarsigma )) = O\Big(   N_{i_0}^{-1} 2^{-q_{i_0}}      \sum_{\bgamma\in 2^{-\bq} \cdot \Gamma^\bot \setminus {\bs}}  
       |\ddot{\cC}_{\bq}(M, \bgamma, i_0)|     
           \Big),
\end{equation}
with
\begin{equation}\nonumber
     \ddot{\cC}_{\bq}(M, \bgamma,i_0)  =  \frac{  (  1-\eta_M(\bgamma))(1-\eta(\gamma_1 /M))  }{\gamma_1}   \prod_{j=2}^s \frac{  
              m ( \gamma_j)  }{\gamma_j}
     \frac{1}{\gamma_{i_0}}  .  
\end{equation}
Applying (\ref{9.10}), we obtain $\max_{\gamma \in \cX, i\in [1,s]} |1/\gamma_i| =O(1)$.

By   (\ref{2.35a})  and (\ref{7.5}), we have
\begin{equation}\nonumber
       \bE(\breve{\cC}_{\bq}(M,\bvarsigma )) =     \bE(\bE_{i_0}(\breve{\cC}_{\bq}(M,\bvarsigma )))
       =O\Big(   N_{i_0}^{-1} 2^{-q_{i_0}}      \sum_{\bgamma\in \cX}  1/|\Nm(\bgamma)|\Big)
        = O(    N_{i_0}^{-1} 2^{-q_{i_0}}   ).
\end{equation}
 Similarly to (\ref{7.23}) - (\ref{7.23a}),   we get  from  (\ref{9.1}) and (\ref{6.4}), that 
\begin{equation}\nonumber
       \bE( \widetilde{\cC}_{10}(M))    =  O\Big(  \sum_{ \substack{\bvarsigma \in \{1,-1\}^s\\  \bvarsigma  \neq -\b1}} 
       \sum_{\bq \in  \cG_{9}}    N_{i_0(\bvarsigma)}^{-1} 2^{-q_{i_0(\bvarsigma)}}  \Big)
\end{equation}
\begin{equation}\nonumber
  =  O\Big(  \sum_{ 1 \leq i \leq s}   \sum_{  j \leq  n+ 0.5\log_2 n} 
       \sum_{\bq \in  \cG_{9}, q_{i} =-j }    2^{-n+j}  \Big)
  = O\Big( n^{s-2}  \sum_{  j \leq   1/2\log_2 n}     2^{j} \Big) =O(n^{s-3/2}).
\end{equation}
Using  (\ref{9.11})  and (\ref{9.14}), we obtain the assertion of  Lemma 19. \qed \\

{\bf Lemma 20.} {\it With notations as above
\begin{equation}\nonumber
     \bE( \widetilde{\cC}_{3}(M))    =     \widetilde{\cC}_{12}(M) +O( n^{s-3/2}), \quad M= [\sqrt{n}],
\end{equation}  
where
\begin{equation}\label{9.19}
     \widetilde{\cC}_{12}(M)  = (-1)^s(2\sqrt{-1})^{-s} \sum_{{\bq} \in \cG_9} \sum_{\bgamma_0 \in \Delta_p}  e(\langle \bgamma_0, \bb/p \rangle)\check{\cC}_{\bq}(\bgamma_0 ), 
\end{equation}
with
\begin{equation} \nonumber
  \check{\cC}_{\bq}(\bgamma_0 ) =    M^{-1} \sum_{\bgamma \in \Gamma_{M,\bq} (\bgamma_0)} g(\bgamma), \quad g(\bx) = 
         \eta(2\Nm (\bx))(1-   \eta(x_1) ))  \MM(\bx)/\Nm(\bx), 
\end{equation}
and
\begin{equation}\nonumber
  \Gamma_{M,\bq}(\bgamma_0)= ( p 2^{-\bq} \cdot \Gamma^{\bot} +\bgamma_0 )\cdot(1/M,1,1,...,1).
\end{equation} }

{\bf Proof.}  By (\ref{9.3a}), (\ref{9.8}) and Lemma 19, it is enough to prove that
\begin{equation}\nonumber
      \widetilde{\cC}_{11}(M)    =     \widetilde{\cC}_{12}(M) +O( n^{s-3/2}).
\end{equation}  
Consider $ \breve{\cC}_{\bq}(M, -\b1)$.
Let
\begin{equation}\nonumber
  \bar{\cC}_{\bq}(M, -\b1) =          \sum_{\bgamma\in \Gamma^\bot \setminus {\bs}}
         (  1-\eta_M(\bgamma))  e(\langle \bgamma, \bb/p \rangle))
\end{equation}
\begin{equation}\nonumber
  \times     \eta(2^{-q_1}\gamma_1/M)     \MM (2^{-\bq} \cdot \bgamma) /\Nm (\bgamma).
\end{equation}
By (\ref{9.5}), we have 
\begin{equation}\nonumber
  |\breve{\cC}_{\bq}(M, -\b1) - \bar{\cC}_{\bq}(M, -\b1)|   \leq           \sum_{\bgamma\in \Gamma^\bot \setminus {\bs}}
         |(  1-\eta_M(\bgamma)) \eta(2^{-q_1}\gamma_1/M) \MM (2^{-\bq} \cdot \bgamma) |
\end{equation}
\begin{equation}\nonumber
    \times   |(\widehat{\Omega}(\tau \bgamma) -1)   /\Nm (\bgamma)|.
\end{equation}
 We examine the case $(1-   \eta(\gamma_1 2^{-q_1 }/M))  \MM(2^{-\bq}  \bgamma) \neq 0$.
By (\ref{2.35a})  and (\ref{6.01a}), we get  $|\gamma_1| \leq M2^{q_1+1}$  and $|\gamma_i| \leq 2^{q_i+2}, \; i\geq 2$.

Hence, we obtain from (\ref{6.4})  and (\ref{9.1}), that  $|\tau \gamma_i| \leq 4n^{-s+1/2}, \; i \geq 1$  for $\bq \in \cG_9$. \\
Applying (\ref{2.13}), we get $\widehat{\Omega}(\tau \bgamma) =1+O(n^{-s+1/2})$ for $\bq \in G_9$.
Bearing in mind (\ref{7.5}), we have
\begin{equation}\label{9.21}
               \breve{\cC}_{\bq}(M, -\b1) =      \bar{\cC}_{\bq}(M, -\b1)  +O(n^{-1}).
\end{equation}
Taking into account that $\eta(0) =0$ (see (\ref{2.35a})), we get
\begin{equation}\nonumber
  \bar{\cC}_{\bq}(M, -\b1) =       \sum_{\bgamma_0 \in \Delta_p}  e(\langle \bgamma_0, \bb/p \rangle)   \acute{\cC}_{\bq}( \bgamma_0),
\end{equation}
with
\begin{equation}\nonumber
    \acute{\cC}_{\bq}( \bgamma_0) =  \sum_{\bgamma\in 2^{-\bq}(p\Gamma^\bot +\bgamma_0)} \eta(2|\Nm(\bgamma)|/M) (1-\eta(\gamma_1/M))     \MM ( \bgamma) /\Nm (\bgamma).
\end{equation}
It is easy to verify that $ \acute{\cC}_{\bq}( \bgamma_0) =  \check{\cC}_{\bq}(\bgamma_0)$.
By (\ref{9.7}) and (\ref{9.19}), we obtain
 \begin{equation}\nonumber
  \widetilde{\cC}_{11}(M) =(-1)^s   (2\sqrt{-1})^{-s}  \sum_{{\bq} \in \cG_9}  
    \Big(  \sum_{\bgamma_0 \in \Delta_p}  e(\langle \bgamma_0, \bb/p \rangle)   \breve{\cC}_{\bq}( \bgamma_0)  +O(n^{-1}) \Big)          = \widetilde{\cC}_{12}(M) + O(n^{s-2}).
\end{equation}  
Hence, Lemma 20 is proved. \qed \\

We consider  Koksma--Hlawka inequality  (see e.g. [DrTy, p. 10, 11]):
 
  \texttt{Definition 5.} {\it
  Let a function $f\; : \; [0,1]^s \to \RR$  have continuous partial derivative $\partial^l f^{(F_l)} /\partial x_{i_1} \cdots  \partial x_{i_l}  $ on 
 on the $s-l$ dimensional face $F_l$, defined by $ x_{i_1} = \cdots = x_{i_l} =1$, and let 
\begin{equation}\nonumber
  V^{(s-l)}(f^{F_l}) = \int_{F_l} \Big| \frac{\partial^l f^{(F_l)}}{\partial x_{i_1} \cdots  \partial x_{i_1}} \Big| d x_{i_1} \cdots d x_{i_l}. 
\end{equation}
Then the number
\begin{equation}\nonumber
  V(f) = \sum_{0 \leq l <s} \sum_{F_l}  V^{(s-l)}(f^{F_l})
\end{equation}
is called a Hardy and Krause variation.} \\

{\bf Theorem F.} (Koksma--Hlawka) {\it  Let $f$ be of bounded variation on
 $[0, 1]^s$ in the sense of Hardy and Krause. 
 Let $((\beta_{k,K})_{k=0}^{K-1})$ be a $K$-point set in an $s$-dimensional unit cube $[0,1)^s$. Then we have}
\begin{equation}\nonumber
    \Big| \frac{1}{K} \sum_{0 \leq k \leq K-1}f(\beta_{k,K}) - \int_{[0,1]^s} f(\bx) d\bx \Big|        \leq   V(f)  D((\beta_{k,K})_{k=0}^{K-1})   .
\end{equation}  \\

{\bf Lemma 21.}  {\it With notations as above} 
\begin{equation}\nonumber
     \bE( \widetilde{\cC}_{3}(M))       =O( n^{s-5/4}),  \quad M= [\sqrt{n}].
\end{equation}  
{\bf Proof.}  By (\ref{9.19}) $g(\bx) =\eta(2\Nm (\bx))(1-   \eta(x_1) ))  \MM(\bx)/\Nm(\bx)$. We have that $g$  is the odd function, with respect to each coordinate, and $g(\bx) =0$ for $\bx \notin [-2,2]\times [-4,4]^{s-1}$. Hence
\begin{equation}\nonumber
 \int_{[-2,2]\times [-4,4]^{s-1} } g(\bx)  d \bx =0.
\end{equation}
Let  $f(\bx) =g((4x_1 -2,8x_2 -4,...,8x_s-4))$. It is easy to verify that $f(\bx) =0$ for $\bx \notin [0, 1]^s$, and
\begin{equation}\nonumber
   \int_{[0, 1]^s} f(\bx)  d \bx= \int_{[-2,2]\times [-4,4]^{s-1} } g(\bx)  d \bx =0.
\end{equation}
We see that $f$ is of bounded variation on
 $[0, 1]^s$ in the sense of Hardy and Krause.
  Let  $ \ddot{\Gamma}(\bgamma_0) = \{ ((\gamma_1 +2)/4,(\gamma_2+4)/8,...,(\gamma_s +4)/8) \; | \;    \bgamma \in \Gamma_{M,\bq} (\bgamma_0) \}$.
 
Using  (\ref{9.19}), we obtain
\begin{equation}\nonumber
  \check{\cC}_{\bq}(\bgamma_0 ) =    M^{-1} \sum_{\bgamma \in \ddot{\Gamma}(\bgamma_0)} f(\bgamma).
\end{equation}
Let 
$H=\ddot{\Gamma}(\bgamma_0)\cap [0,1)^{s}$, and  $K =\#H$. Applying Theorem A, we get $K \in [c_1M,c_2M]$ for some $c_1,c_2 >0$.
We enumerate the set $H$ by a sequence $((\beta_{k,K})_{k=0}^{K-1})$.

 By Theorem A, we have  $D((\beta_{k,K})_{k=0}^{K-1}) =O(M^{-1} \ln^{s-1} M)$. 
 
 Using  Theorem F, we obtain $  \check{\cC}_{\bq}(\bgamma_0 ) = O(M^{-1} \ln^{s-1} M) $.  

 Bearing in mind that $ \# \cG_3 =O(n^{s-1})$ (see (\ref{7.17a})), we derive from (\ref{9.19}) that 
$ \widetilde{\cC}_{12}(M) =O(n^{s-1} M^{-1} \ln^{s-1} M )$.

Applying Lemma 20, we obtain the assertion of the Lemma 21. \qed \\

Now using (\ref{6.28}), Corollary 1 and  Lemma 21, we get  \\

{\bf  Corollary 2.} {\it With notations as above} 
\begin{equation}\nonumber
       \bE(\cB(\bb/p,M)) =  O( n^{s-5/4}), \quad M= [\sqrt{n}].
\end{equation}\\

Let $\bN = (N_1,...,N_s)$,  $N = N_1 \cdots N_s$, $n=s^{-1}\log_2 N$, $c_9=   0.25 ( \pi^s  \det \Gamma   )^{-1} c_{8}$ and $M = [\sqrt{n}]$. From  Lemma 12, Corollary 2 and (\ref{2.46}),   we obtain that there exist $N_0 >0$, and   $\bb \in \Delta_p$ such that 
\begin{equation}\label{9.25}
    \sup_{\btheta \in [0,1]^s} | \bE( \cR) (B_{\btheta \cdot\bN}+\bb/p,\Gamma )  |  \geq
         c_{9}n^{s-1}    \quad {\rm for} \quad N >N_0.
\end{equation} \\

{\bf 2.10. End of proof.} 
{\it End of the proof of Theorem 1.}\\
We set $ \widetilde{\cR} (\bz,\by) =\cR (B_{\by-\bz}+\bz,\Gamma)$, where $y_i \geq z_i$ $(i=1,...,s)$  (see (\ref{1.2})).
Let us introduce the difference operator $ \dot{\Delta}_{a_i,h_i} \; : \; \RR^s  \to \RR $, defined by the
formula
\begin{equation}\nonumber
  \dot{\Delta}_{a_i,h_i} \tilde{\cR}(\bz,\by)  = \tilde{\cR}(\bz,(y_1,...,y_{i-1},h_i,y_{i+1},...,y_s))-\tilde{\cR}(\bz,(y_1,...,y_{i-1},a_i,y_{i+1},...,y_s)).
\end{equation}
Similarly to [Sh, p. 160,  ref.7], we derive
\begin{equation} \label{10.1}
  \dot{\Delta}_{a_1,h_1}   \cdots \dot{\Delta}_{a_s,h_s} \widetilde{\cR}(\bz,\by)  = \widetilde{\cR}(\ba,{\bf h}),
\end{equation}
where $h_i \geq a_i \geq z_i $ $(i=1,...,s)$.
Let $\bff_1,...,\bff_s$ be a basis of $\Gamma$. We have that $F= \{ \rho_1 \bff_1+ \cdots +  \rho_s\bff_s \; | \; (\rho_1,...,\rho_s) \in [0,1)^s     \}$ is  the fundamental set of $\Gamma$.
It is easy to see that $ \cR (B_{\bN}+\bx,\Gamma) = \cR (B_{\bN}+\bx+\bgamma,\Gamma)$ for all $\bgamma \in \Gamma$.
Hence, we can assume  in Theorem 1 that $\bx \in F$. Similarly, we can assume  in Corollary 2 that $\bb/p \in F$.    We get that there exists $\bgamma_0 \in  \Gamma$ with $|\bgamma_0| \leq 4 \max_i |\bff_i|$ and $x_i < (\bb/p)_i + \gamma_{0,i}$, $i=1,...,s$.
 Let $\bb_1= \bb+p\bgamma_0$. 
By (\ref{9.25}), we have that  there exists $\btheta \in [0,1]^s$ and $\bb \in \Delta_p$ such that
\begin{equation}\label{10.2}
   |\widetilde{\cR}(\bb_1/p,  \bb_1/p + \btheta \cdot \bN )|  \geq c_{9}  n^{s-1}.
\end{equation}
Let $\cS =\{ \by \; | \; y_i =(\bb_/p)_i, (\bb_/p)_i +\theta_i  N_i, \; i=1,...,s \}$. We see $\# \cS =2^s$. 
From  (\ref{10.1}),  we obtain that $\widetilde{\cR}(\bb_1/p,  \bb_1/p + \btheta \cdot \bN )$ is the sum of $2^s$ numbers
 $\pm  \widetilde{\cR}(\bx,  \by^{j})$, where $\by^{j} \in \cS$. By (\ref{10.2}), we get   
\begin{equation}\nonumber
    |\cR (B_{\by-\bx}+\bx,\Gamma)| = |\widetilde{\cR}(\bx,  \by)|  \geq 2^{-s}c_{9} n^{s-1}, \quad  {\rm for \; some} \quad \by \in \cS .
\end{equation} 
Therefore, Theorem 1 is proved. \qed  \\

{\it  Proof of Theorem 2.}
Let $n \geq 1,\; N \in [2^n, 2^{n+1})$,  $\by =(y_1,...,y_{s-1})$ and $\Gamma=\Gamma_{\cM}$.
By (\ref{1.14}) and [Le3, p.41], we have
\begin{equation}\label{9.40}
  ([y_s N] +1) \Delta(B_{\by}, (\beta_{k,N})_{k=0}^{[y_s N]}  ) = \alpha_1 -   y_1  \cdots y_{s-1} \alpha_2 +O(\log_2^{s-1}n),
\end{equation}
where
\begin{equation}\nonumber
 \alpha_1 =  \cN (B_{(y_1,...,y_{s-1},y_sN  \det \Gamma)}+\bx,\Gamma), \;\;
 {\rm and} \;\;  \alpha_2 =  
  \cN (B_{(1,...,1,y_sN \det \Gamma)}+\bx,\Gamma).
\end{equation}
From (\ref{1.2}), we get
\begin{equation}\label{9.42}
  \alpha_1 - y_1  \cdots y_{s-1}\alpha_2 =  \beta_1 -  y_1  \cdots y_{s-1}\beta_2, 
\end{equation}
with
\begin{equation}\nonumber
 \beta_1 =  \cR (B_{(y_1,...,y_{s-1},y_sN  \det \Gamma)}+\bx,\Gamma), \;\;
 {\rm and} \;\;  \beta_2 =  
   \cR (B_{(1,...,1,y_sN \det \Gamma)}+\bx,\Gamma).
\end{equation}
Let $y_0 =0.125\min\big(1, 1/\det \Gamma, (c_1(\cM)/c_0(\Gamma))^{1/(s-1)} \big)$, $\btheta =(\theta_1,...,\theta_{s})$,  $y_i = y_0 \theta_i$,  $i=1,...,s-1$, and $y_s=\theta_s$.  Using Theorem A, we obtain 
\begin{equation}\nonumber
    |y_1\cdots y_{s-1} \cR (B_{(1,...,1,y_sN \det \Gamma)}+\bx,\Gamma)| \leq y_0^{s-1} c_0(\Gamma) \log_2^{s-1}(2+y_sN \det \Gamma)
\end{equation}
\begin{equation}\label{9.44}
    \leq (2y_0)^{s-1} c_0(\Gamma) \log_2^{s-1}N \leq  0.25c_1(\cM) n^{s-1}  \quad  {\rm for} \quad  N > \det \Gamma +2.
\end{equation} 
Applying Theorem 1, we have 
\begin{equation}\nonumber
  \sup_{\btheta \in [0,1)^{s}} |\cR (B_{(\theta_1 y_0 ,...,\theta_{s-1}y_0, \theta_s N  \det \Gamma)}+\bx,\Gamma)| \geq c_1 (\cM)\log_2^{s-1}(y_0^{s-1} \det \Gamma N)
\end{equation}
\begin{equation}\nonumber
 \geq c_1 (\cM) n^{s-1}(1 +n^{-1}(s-1)\log_2 (y_0^{s-1} \det \Gamma )) \geq  0.5c_1(\cM) n^{s-1} 
\end{equation}
for $n> 10(s-1)|\log_2 (y_0^{s-1} \det \Gamma )|$. 
Using (\ref{1.16}), (\ref{9.40}), (\ref{9.42}) and  (\ref{9.44}), we get the assertion of Theorem 2. \qed  \\

{\bf Bibliography.}

 [Be] Beck, J., A two-dimensional van Aardenne-Ehrenfest theorem in
irregularities of distribution, Compos. Math. 72 (1989), no. 3, 269-339.

[BC] Beck, J., Chen, W.W.L.
 Irregularities of Distribution,  
 Cambridge Univ. Press,  Cambridge 1987.

[Bi]  Bilyk, D.,  On Roth's orthogonal function method in discrepancy theory. Unif. Distrib. Theory 6 (2011), no. 1, 143-184. 

[BiLa] Bilyk, D., and  Lacey, M.,   The Supremum Norm of the Discrepancy
Function: Recent Results and Connections, 	arXiv:1207.6659

[BS]  Borevich, A. I., Shafarevich, I. R., Number Theory, Academic Press, New York, 1966. 

[ChTr]  Chen, W., Travaglini G., Some of Roth's Ideas in discrepancy theory. Analytic Number Theory, 150-163, Cambridge Univ. Press, Cambridge, 2009. 

[Ch] Chevalley, C., Deux th\'eor\`{e}mes d'arithm\'etique, J. Math. Soc. Japan {\bf 3} (1951) 36-44.

[DrTi]  Drmota, M.,  Tichy, R., Sequences, Discrepancies and Applications, Lecture Notes in Mathematics
1651, 1997.

 [GL]  Gruber, P.M,  Lekkerkerker, C.G., Geometry of Numbers, North-Holland, New-York,   1987.

 [Is] Isaacs, I. M., Character Theory of Finite Groups,  AMS Chelsea Publishing, Providence, RI, 2006.

 [KaNi]  Katok, A., Nitica, V.,  Rigidity in Higher Rank Abelian Group Actions, Volume I, Introduction and Cocycle Problem,  Cambridge University Press, Cambridge, 2011.

[KO]  Kostrikin, A. I. Introduction to Algebra, Springer-Verlag, New York-Berlin, 1982.

[La1] Lang, S., Algebraic Number Theory, Springer-Verlag, New York, 1994. 

 [La2] Lang, S., Algebra, Springer-Verlag, New York, 2002.

 [Le1] Levin, M.B., On low discrepancy sequences and low
discrepancy ergodic transformations of the multidimensional unit
cube,    Israel J.  Math. {\bf 178} (2010),  61-106. 

[Le2]   Levin, M.B., Adelic constructions of low discrepancy sequences, Online J. Anal. Comb. No. 5 (2010), 27 pp. 

 [Le3]  Levin, M.B., On Gaussian limiting  distribution of lattice points in 
  a parallelepiped, Arxiv 
   
[MuEs]  Murty, M.R., Esmonde, J., Problems in algebraic number theory,  Springer-Verlag, New York, 2005.
 
[Na]  Narkiewicz, W., Elementary and Analytic Theory of Algebraic Numbers,  Springer-Verlag, Berlin, 1990.

[Ne] Neukirch, J. Algebraic Number Theory,  Springer-Verlag, Berlin, 1999.
  
[NiSkr]  Nikishin, N. A.; Skriganov, M. M.,  On the distribution of algebraic numbers in parallelotopes. (Russian) Algebra i Analiz 10 (1998), no. 1, 68--87; translation in St. Petersburg Math. J. 10 (1999), no. 1, 53-68 
  
 [Sh] Shiryaev, A.N., Probability, Springer-Verlag, New York, 1996.
  
 [Skr] Skriganov, M.M.,  Construction of uniform distributions in terms of
geometry of numbers,
 Algebra i Analiz {\bf6}, no.~3 (1994), 200-230; 
Reprinted in St.~Petersburg Math. J. {\bf6}, no.~3 (1995), 635-664.

   [SW]  Stein, E., Weiss, G.,  Introduction to Fourier Analysis on Euclidean Spaces, Princeton University Press, New-York, 1971.
   
[Ve]    Veech, W.A., Periodic points and invariant pseudomeasures for toral endomorphisms, Ergodic Theory Dynam. Systems 6 (1986), no. 3, 449-473. 

[Wi] Wills, J. M., Zur Gitterpunktanzahl konvexer Mengen, 
Elem. Math.  28 (1973), 57-63.
\\
\\
{\bf Address}: Department of Mathematics,
Bar-Ilan University, Ramat-Gan, 5290002, Israel \\
{\bf E-mail}: mlevin@math.biu.ac.il\\
\end{document}